# TRANSPORTATION COST-INFORMATION INEQUALITIES AND APPLICATIONS TO RANDOM DYNAMICAL SYSTEMS AND DIFFUSIONS

By H. Djellout, A. Guillin and L. Wu

*Université Blaise Pascal, Université Paris IX Dauphine and Université Blaise Pascal*

We first give a characterization of the $L^1$-transportation cost-information inequality on a metric space and next find some appropriate sufficient condition to transportation cost-information inequalities for dependent sequences. Applications to random dynamical systems and diffusions are studied.

**1. Introduction and questions.** Let $(E, d)$ be a metric space equipped with $\sigma$-field $\mathcal{B}$ such that $d(\cdot, \cdot)$ is $\mathcal{B} \times \mathcal{B}$-measurable. Given $p \geq 1$ and two probability measures $\mu$ and $\nu$ on $E$, we define the quantity

$$(1.1) \qquad W_p^d(\mu, \nu) = \inf \left( \iint d(x, y)^p \, d\pi(x, y) \right)^{1/p},$$

where the infimum is taken over all probability measures $\pi$ on the product space $E \times E$ with marginal distributions $\mu$ and $\nu$ [say coupling of $(\mu, \nu)$]. This infimum is finite as soon as $\mu$ and $\nu$ have finite moments of order $p$. This quantity is commonly referred to as $L^p$-Wasserstein distance between $\mu$ and $\nu$. When $d$ is the trivial metric $(d(x, y) = \mathbb{1}_{x \neq y})$, $2W_1^d(\mu, \nu) = \|\mu - \nu\|_{\text{TV}}$, the total variation of $\mu - \nu$.

The Kullback information (or relative entropy) of $\nu$ with respect to $\mu$ is defined as

$$(1.2) \qquad H(\nu/\mu) = \begin{cases} \int \log \dfrac{d\nu}{d\mu} \, d\nu, & \text{if } \nu \ll \mu, \\ +\infty, & \text{otherwise.} \end{cases}$$

We say that the probability measure $\mu$ satisfies the $L^p$-transportation cost-information inequality on $(E, d)$ if there is some constant $C > 0$ such that for any probability measure $\nu$,

$$(1.3) \qquad W_p^d(\mu, \nu) \leq \sqrt{2CH(\nu/\mu)}.$$

To be short, we write $\mu \in T_p(C)$ for this relation.

The cases "$p = 1$" and "$p = 2$" are particularly interesting. That $T_1(C)$ are related to the phenomenon of measure concentration was emphasized by Marton [10, 11], Talagrand [18], Bobkov and Götze [2] and amply explored by Ledoux [8, 9]. The $T_2(C)$, first established by Talagrand [18] for the Gaussian measure, has been brought into relation with the log-Sobolev inequality, Poincaré inequality, inf-convolution, Hamilton–Jacobi's equations by Otto and Villani [15] and Bobkov, Gentil and Ledoux [1]. Since those important works, a main trend in the field is to put on relations of $T_p(C)$ with other functional inequalities (of geometrical nature in particular). In this paper we shall study three questions around the following problem going somehow to the opposite direction: *how to establish the "$T_p(C)$" without reference to other functional inequalities in various concrete situations?*

To raise our first question, let us mention the following:

THEOREM 1.1 (Bobkov and Götze [2]). *$\mu$ satisfies the $L^1$-transportation cost-information inequality on $(E, d)$ with constant $C > 0$, that is, $\mu \in T_1(C)$, if and only if for any Lipschitzian function $F : (E, d) \to \mathbb{R}$, $F$ is $\mu$-integrable and*

$$(1.4) \qquad \int_E e^{\lambda(F - \langle F \rangle_\mu)} \, d\mu \leq \exp\left(\frac{\lambda^2}{2} C \|F\|_{\mathrm{Lip}}^2\right) \qquad \forall \lambda \in \mathbb{R},$$

*where $\langle F \rangle_\mu = \int_E F \, d\mu$ and*

$$\|F\|_{\mathrm{Lip}} = \sup_{x \neq y} \frac{|F(x) - F(y)|}{d(x, y)} < +\infty.$$

*In that case,*

$$\mu(F - \langle F \rangle_\mu > r) \leq \exp\left(-\frac{r^2}{2C\|F\|_{\mathrm{Lip}}^2}\right) \qquad \forall r > 0.$$

It might be worthwhile to recall the classical Pinsker–Csizsar inequality which is the starting point of many recent works. By the coupling characterization of the total variation distance $\|\cdot\|_{\mathrm{TV}}$, the Pinsker–Csizsar inequality

$$\|\nu - \mu\|_{\mathrm{TV}} \leq \sqrt{\tfrac{1}{2} H(\nu/\mu)}$$

says that w.r.t. the trivial distance $d(x, y) = \mathbb{1}_{x \neq y}$ on $E$, any probability measure $\mu$ satisfies the $L^1$-transportation cost-information inequality with the sharp constant $C = 1/4$. And by Theorem 1.1, the Pinsker–Csiszar inequality for the trivial distance follows from the classical well-known inequality: for a real bounded random variable $\xi$ with values in $[a, b]$,

$$\mathbb{E} e^{\xi - \mathbb{E}\xi} \leq \exp\left(\frac{(b-a)^2}{8}\right)$$

(and vice versa).

We now do a simple remark. Assume that $\mu \in T_1(C)$ or, equivalently, (1.4). Let $\gamma(d\lambda)$ be the standard Gaussian law $\mathcal{N}(0, 1)$ on $\mathbb{R}$. We have for any Lipschitzian function $F$ on $E$ with $\langle F \rangle_\mu = 0$ and $\|F\|_{\text{Lip}} \leq 1$, and $a \in \mathbb{R}$,

$$\int_E \exp\left(\frac{a^2}{2} F^2\right) d\mu = \int_E \int_\mathbb{R} e^{a\lambda F} \gamma(d\lambda) \, d\mu \leq \int_\mathbb{R} \exp\left(\frac{C}{2} a^2 \lambda^2\right) \gamma(d\lambda)$$

$$= \begin{cases} \dfrac{1}{\sqrt{1 - a^2 C}}, & \text{if } \dfrac{a^2}{2} < \dfrac{1}{2C}, \\ +\infty, & \text{otherwise.} \end{cases}$$

Applying it to $F(x) := d(x, x_0) - \int d(x, x_0) \, d\mu(x)$, we obtain

$$\int e^{cd^2(x, x_0)} \, d\mu(x) < +\infty \qquad \forall c \in \left(0, \frac{1}{2C}\right).$$

In particular, for all $\delta \in (0, \frac{1}{4C})$ we have,

$$(1.5) \qquad\qquad \iint e^{\delta \, d^2(x, y)} \, d\mu(x) \, d\mu(y) < +\infty.$$

That naturally leads to the following questions:

QUESTION 1. Will the Gaussian tail (1.5) be sufficient for the $L^1$-transportation cost-information inequality of $\mu$?

The second question is about dependent tensorizations of the $T_p(C)$. Let, for example, $\mathbb{P}_x^n$, the law of a homogeneous Markov chain $(X_k(x))_{1 \leq k \leq n}$ on $E^n$ starting from $x \in E$, with transition kernel $P(x, dy)$.

QUESTION 2. Assume that $P(x, \cdot) \in T_p(C)$ for all $x \in E$. Where is the appropriate condition under which $\mathbb{P}_x^n$ satisfies the $L^p$-transportation cost-information inequality w.r.t. the metric

$$d_{l_p}(x, y) := \left(\sum_{i=1}^n d(x_i, y_i)^p\right)^{1/p} ?$$

The same question can be raised for the law of an arbitrary dependent sequence $(X_k)_{1 \leq k \leq n}$. When $(X_k)_{1 \leq k \leq n}$ are independent, this question has a rapid and affirmative answer, see [8, 9] and references therein.

In the dependent case, when $d$ is the trivial metric, and $p = 1$ (and $d_{l_1}$ becomes the Hamming distance on $E^n$), Marton [10] generalized the Pinsker–Csiszar inequality to the law of the so called "contracting" Markov chains:

$$(1.6) \qquad \frac{1}{2} \sup_{(x_{i-1}, y_{i-1})} \|P_i(\cdot/y_{i-1}) - P_i(\cdot/x_{i-1})\|_{\mathrm{TV}} \leq r < 1.$$

Her approach is based on coupling ideas, natural by the definition of the involved Wasserstein distance. Her results have been strengthened by Marton [11, 12] and Dembo [4] and have been generalized to uniform mixing processes by Samson [17] and Rio [16].

However, the trivial distance does not reflect the natural metric structure of the state space $E$ to which usual Markov processes such as random dynamical systems or diffusions are related and that is why the uniform mixing assumption was made in her work (and also in [17]). This is a main motivation for Question 2.

For the $L^2$-transportation cost-information inequality $T_2(C)$, recall that Talagrand [18] proved that the standard Gaussian law $\gamma = \mathcal{N}(0, 1)$ satisfies $T_2(C)$ on $\mathbb{R}$ w.r.t. the Euclidean distance with the sharp constant $C = 1$ and found that $T_2(C)$ is stable for product (or independent) tensorization. To our knowledge the Markovian tensorization of $T_2(C)$ has not been investigated in the literature.

Since the works of Otto and Villani [15] and Bobkov, Gentil and Ledoux [1], we know that $T_2(C)$ follows from the log-Sobolev inequality in the framework of Riemannian manifolds. Indeed, all known $T_2(C)$-inequalities up to now can be derived from the log-Sobolev inequality. An important open question in the field is whether $T_2(C)$ is strictly weaker than the log-Sobolev inequality. Hence, it would be interesting to investigate the following question:

QUESTION 3. How do we establish the $T_2(C)$-inequality in situations where the log-Sobolev inequality is unknown?

This paper is written around those three questions and it is organized as follows. The next section is the general theoretical part of this paper. After noticing the stability of $T_p(C)$ under Lipschitzian map and under weak convergence in Sections 2.1 and 2.2, in Section 2.3 we prove that condition (1.5) is, in fact, sufficient for the $L^1$-transportation cost-information inequality, solving Question 1. In Section 2.4 we revisit the coupling method of Marton and show that it actually works for dependent tensorization of $T_p(C)$ for $1 \leq p \leq 2$, under a contraction assumption [see (C1) in Theorem 2.5] close to Marton's (1.6). Section 2.5 is devoted to revisit the McDiarmid–Rio martingale method which allows us to obtain a much more subtle condition (C1$'$) than (C1) for tensorization of $T_1(C)$ in Theorem 2.11.

Sections 3 and 4 contain several applications of the general results in Section 2 to random dynamical systems and diffusions which are our main motivation for Question 2.

In Section 5, quite independent, we present a direct approach of $T_2(C)$ for diffusions, by means of the Girsanov transformation, with respect to the usual Cameron–Martin metric or $L^2$-metric.

The reader may consult the recent monograph by Villani [19] for an extended (analytical and geometrical) treatment on transportation.

**2. Criteria for $T_p(C)$.** Throughout this paper let $(E, d)$ be a metric space equipped with $\sigma$-field $\mathcal{B}$ such that $d(\cdot, \cdot)$ is $\mathcal{B} \times \mathcal{B}$-measurable; and when $(E, d)$ is separable, $\mathcal{B}$ will be the Borel $\sigma$-field.

2.1. *Stability under push-forward by Lipschitz map.* We begin with the stability of $T_p(C)$ under Lipschitzian map and under weak convergence, which will be useful later.

LEMMA 2.1. *Assume that $\mu \in T_p(C)$ on $(E, d_E)$ and $(F, d_F)$ is another metric space. If $\Psi : (E, d_E) \to (F, d_F)$ is Lipschitzian,*

$$d_F(\Psi(x), \Psi(y)) \leq \alpha \, d_E(x, y) \qquad \forall x, y \in E,$$

*then $\tilde{\mu} := \mu \circ \Psi^{-1} \in T_p(C\alpha^2)$ on $(F, d_F)$.*

PROOF. Let $\tilde{\nu}$ be a probability measure such that $H(\tilde{\nu}/\tilde{\mu}) < +\infty$. The key remark is

(2.1) $$H(\tilde{\nu}/\tilde{\mu}) = \inf\{H(\nu/\mu); \nu \circ \Psi^{-1} = \tilde{\nu}\}.$$

To prove it, putting $\nu_0(dx) := \frac{d\tilde{\nu}}{d\tilde{\mu}}(\Psi(x))\mu(dx)$, we see that $\nu_0 \circ \Psi^{-1} = \tilde{\nu}$. We have for any $\nu$ so that $\nu \circ \Psi^{-1} = \tilde{\nu}$,

$$H(\nu/\mu) = H(\nu_0/\mu) + \int d\tilde{\nu}(y) H(\nu_y/\mu_y),$$

where $\nu_y := \nu(\cdot/\Psi = y)$ and $\mu_y := \mu(\cdot/\Psi = y)$ are, respectively, the regular conditional distribution of $\nu$, $\mu$ knowing $\Psi = y$. Hence, (2.1) follows.

With (2.1) in hand, the rest of the proof is easy and is omitted. □

2.2. *Stability under weak convergence.*

LEMMA 2.2. *Let $(E, d)$ be a metric, separable and complete space (Polish, say) and $(\mu_n, \mu)_{n \in \mathbb{N}}$ a family of probability measures on $E$. Assume that $\mu_n \in T_p(C)$ for all $n \in \mathbb{N}$ and $\mu_n \to \mu$ weakly. Then $\mu \in T_p(C)$.*

PROOF. Recall at first two facts (see, e.g., [19]):

1. If $\mu_n \to \mu$ and $\nu_n \to \nu$ weakly, then $\liminf_{n \to \infty} W_p(\mu_n, \nu_n) \geq W_p(\mu, \nu)$.
2. If $\mu_n \to \mu$ weakly and $\{d(x, x_0)^p, \mu_n(dx)\}$ is uniformly integrable, $W_p(\mu_n, \mu) \to 0$.

What one needs to prove is

$$W_p^2(f\mu, \mu) \leq 2C \int f \log f \, d\mu$$

for all $f$ such that $f\mu$ is a probability. By approximation (and using fact 2 above), it is sufficient to prove the result for continuous $f$ so that $1/N \leq f \leq N$ over $E$ for some $N \geq 1$. Let $a_n = \int f \, d\mu_n$ and we have by "$\mu_n \in T_p(C)$,"

$$W_p^2\left(\frac{f\mu_n}{a_n}, \mu_n\right) \leq 2C \int \left(\frac{f}{a_n}\right) \log\left(\frac{f}{a_n}\right) d\mu_n = \frac{2C}{a_n} \int f \log f \, d\mu_n.$$

Since $\mu_n$ converges weakly to $\mu$, $a_n$ converges to $\mu(f) = 1$, and one can pass to the limit in the right-hand side of this last inequality. For the convergence of the left-hand side, it is enough to apply the lower semi-continuity of $W_p$. □

2.3. *Characterization of $T_1(C)$ by "Gaussian tail."* We present here a characterization of $T_1(C)$, based on the Bobkov and Götze [2] result, that is, some Gaussian integrability property.

THEOREM 2.3. *A given probability measure $\mu$ on $(E, d)$ satisfies the $L^1$-transportation cost-information inequality with some constant $C$ on $(E, d)$ if and only if (1.5) holds. In the latter case,*

$$(2.2) \quad C \leq \frac{2}{\delta} \sup_{k \geq 1}\left(\frac{(k!)^2}{(2k)!}\right)^{1/k} \cdot \left[\iint e^{\delta d^2(x,y)} \, d\mu(x) \, d\mu(y)\right]^{1/k} < +\infty.$$

PROOF. It is enough to show the sufficiency. By Bobkov–Götze's Theorem 1.1, it is enough to show that there is some constant $C = C(\delta)$ verifying (2.2) such that

$$(2.3) \quad \mathbb{E}e^{\lambda F(\xi)} \leq \exp\left(\frac{C\lambda^2}{2}\right) \qquad \forall \lambda \in \mathbb{R},$$

for all $F: E \to \mathbb{R}$ with $\|F\|_{\mathrm{Lip}} \leq 1$ and $\mathbb{E}F(\xi) = 0$, where $\xi$ is a random variable valued in $E$ with law $\mu$, defined on some probability space $(\Omega, \mathcal{F}, \mathbb{P})$.

Let $\xi'$ be an independent copy of $\xi$, defined on the same probability space $(\Omega, \mathcal{F}, \mathbb{P})$. Since $\mathbb{E}F(\xi') = 0$, by the convexity of the $x \to e^x$, we have $\mathbb{E}(e^{-\lambda F(\xi')}) \geq 1$. Consequently, noting that $\mathbb{E}[F(\xi) - F(\xi')]^{2k+1} = 0$, we have

$$\mathbb{E}(e^{\lambda F(\xi)}) \leq \mathbb{E}(e^{\lambda F(\xi)})\mathbb{E}(e^{-\lambda F(\xi')})$$
$$= \mathbb{E}e^{\lambda(F(\xi) - F(\xi'))}$$

$$= 1 + \sum_{k=1}^{\infty} \frac{\lambda^{2k} \mathbb{E}[F(\xi) - F(\xi')]^{2k}}{(2k)!}$$

$$\leq 1 + \sum_{k=1}^{\infty} \frac{\lambda^{2k} \mathbb{E}\, d(\xi, \xi')^{2k}}{(2k)!}.$$

Hence, putting

$$C := 2 \sup_{k \geq 1} \left( \frac{k! \cdot \mathbb{E}\, d(\xi, \xi')^{2k}}{(2k)!} \right)^{1/k},$$

we get

$$\mathbb{E}(e^{\lambda F(\xi)}) \leq 1 + \sum_{k=1}^{\infty} \frac{\lambda^{2k}}{k!} \cdot \left( \frac{C}{2} \right)^k = \exp \left( \frac{C\lambda^2}{2} \right).$$

Thus, for (2.3), it remains to estimate $C$ defined above. Since

$$\mathbb{E} d(\xi, \xi')^{2k} \leq k! \cdot \left( \frac{1}{\delta} \right)^k \mathbb{E} \exp(\delta d(\xi, \xi')^2),$$

we get

$$C \leq \frac{2}{\delta} \sup_{k \geq 1} \left( \frac{(k!)^2}{(2k)!} \right)^{1/k} \cdot \left[ \mathbb{E} \exp(\delta(d(\xi, \xi')^2)) \right]^{1/k} < +\infty$$

the desired estimate (2.2). □

REMARK 2.4. For comparison notice that the Bernoulli measure $\mu$ on $\{0, 1\}$ with $\mu(1) \in (0, 1)$ satisfies $T_1(1/4)$ w.r.t. the trivial metric, but does not satisfy $T_p(C)$ for any $p > 1$ (see [7]). Hence, any probability measure $\mu$ which is not a Dirac measure on $E$ does not satisfy $T_p(C)$ for any $p > 1$ w.r.t. the trivial metric. Another example for illustrating difference of $T_1$ and $T_2$ inequalities is the following.

Let $\mu = \phi(x)^2\, dx$ on $\mathbb{R}$ with $0 \leq \phi \in C_0^{\infty}(\mathbb{R})$ (compact support). It satisfies always $T_1(C)$ w.r.t. the Euclidean $d(x, y) := |y - x|$ by the theorem above. But if the support of $\mu$ (or of $\phi$) has two connected components $I_1, I_2$ with $\mathrm{dist}(I_1, I_2) > 0$, then the corresponding $T_2(C)$ fails. In fact, if contrary to $\mu \in T_2(C)$, then by [15] or [1] the following Poincaré inequality holds:

$$\mathrm{Var}_{\mu}(f) \leq C \int_{\mathbb{R}} f'^2\, d\mu \qquad \forall f \in C_0^{\infty}(\mathbb{R}).$$

Choose now $f$ smooth enough and equal to 1 on $I_1$ and 0 on $I_2$. Then the right-hand side in the Poincaré inequality is 0, whereas the variance of $f$ will be non zero so that the Poincaré inequality cannot hold, neither $T_2(C)$.

This example shows, moreover, that $T_1(C)$ on $\mathbb{R}$ does not imply the Poincaré inequality, unlike $T_2(C)$.

The next two sections are dedicated to the tensorization of $T_p(C)$ for dependent sequences.

2.4. *Weakly dependent sequences: Marton's coupling revisited.* Let $\mathbb{P}$ be a probability measure on the product space $(E^n, \mathcal{B}^n)$, $n \geq 2$. For any $x \in E^n$, $x^i := (x_1, \ldots, x_i)$. Let $P_i(\cdot/x^{i-1})$ denote the regular conditional law of $x_i$ given $x^{i-1}$ for $i \geq 2$ (assume its existence). By convention $P_1(\cdot/x^0)$ is the law of $x_1$ under $\mathbb{P}$, where $x^0 = x_0$ is some fixed point. When $\mathbb{P}$ is Markov, then $P_i(\cdot/x^{i-1}) = P_i(\cdot/x_{i-1})$ is the transition kernel at step $i-1$.

Our aim in this section is to extend transportation cost-information inequalities (1.3) for a probability measure $\mathbb{P}$ on $(E^n, d_{l_p})$, where

$$d_{l_p}(x,y) := \left( \sum_{i=1}^n d(x_i, y_i)^p \right)^{1/p}.$$

THEOREM 2.5.   *Let $\mathbb{P}$ be a probability measure on $E^n$, and $1 \leq p \leq 2$. Assume that $P_i(\cdot/x^{i-1}) \in T_p(C)$ on $(E, d)$ for all $i \geq 1$, $x^{i-1}$ in $E^{i-1}$ ($E^0 := \{x_0\}$). If*

(C1)  *there exist $a_j \geq 0$ with $r^p := \sum_{j=1}^\infty (a_j)^p < 1$ such that*

$$(2.4) \qquad [W_p^d(P_i(\cdot/x^{i-1}), P_i(\cdot/\tilde{x}^{i-1}))]^p \leq \sum_{j=1}^{i-1} (a_j)^p d^p(x_{i-j}, \tilde{x}_{i-j}),$$

*for all $i \geq 1$, $x^{i-1}, \tilde{x}^{i-1}$ in $E^{i-1}$, then for any probability measure $\mathbb{Q}$ on $E^n$,*

$$W_p^{d_{l_p}}(\mathbb{Q}, \mathbb{P}) \leq \frac{1}{1-r} \sqrt{2Cn^{2/p-1} H(\mathbb{Q}/\mathbb{P})}.$$

PROOF.   The proof is similar to the one used for the Hamming distance by Marton [10], however, we have to use the assumption $P_i(\cdot/x^{i-1}) \in T_p(C)$ instead of Pinsker's inequality. Assume that $H(\mathbb{Q}/\mathbb{P}) < +\infty$ (trivial otherwise).

Let $Q_i(\cdot/x^{i-1})$ be the regular conditional law of $x_i$ knowing $x^{i-1}$ for $i \geq 2$ and $Q_1(\cdot/x^0)$ the law of $x_1$, both under law $\mathbb{Q}$. We shall use the Kullback information between conditional distributions,

$$H_i(\tilde{x}^{i-1}) = H(Q_i(\cdot/\tilde{x}^{i-1})/P_i(\cdot/\tilde{x}^{i-1})),$$

and exploit the following important identity:

$$(2.5) \qquad H(\mathbb{Q}/\mathbb{P}) = \sum_{i=1}^n \int_{E^n} H_i(\tilde{x}^{i-1}) \, d\mathbb{Q}(\tilde{x}).$$

The key is to construct an appropriate coupling of $\mathbb{Q}$ and $\mathbb{P}$, that is, two random sequences $\tilde{X}^n$ and $X^n$ distributed according to $\mathbb{Q}$ and $\mathbb{P}$, respectively, on some probability space $(\Omega, \mathcal{F}, \mathrm{P})$.

We define a joint distribution $\mathcal{L}(\tilde{X}^n, X^n)$ by induction as follows. Add artificially time 0 and put $X_0 = \tilde{X}_0 = \tilde{x}^0 = x^0$, the fixed point. Assume that for some $i$, $1 \le i \le n$, $\mathcal{L}(\tilde{X}^{i-1}, X^{i-1})$ is already defined. We have to define the joint conditional distribution $\mathcal{L}(\tilde{X}_i, X_i / \tilde{X}^{i-1} = \tilde{x}^{i-1}, X^{i-1} = x^{i-1})$, where $(\tilde{x}^{i-1}, x^{i-1})$ is fixed (but arbitrary).

Given $\varepsilon > 0$ so small that $r(1+\varepsilon) < 1$, this distribution will have marginal laws

$$\mathcal{L}(\tilde{X}_i / \tilde{X}^{i-1} = \tilde{x}^{i-1}, X^{i-1} = x^{i-1}) = Q_i(\cdot / \tilde{x}^{i-1})$$

and

$$\mathcal{L}(X_i / \tilde{X}^{i-1} = \tilde{x}^{i-1}, X^{i-1} = x^{i-1}) = P_i(\cdot / x^{i-1})$$

so as to satisfy

$$\mathbb{E}(d(\tilde{X}_i, X_i)^p / \tilde{X}^{i-1} = \tilde{x}^{i-1}, X^{i-1} = x^{i-1})$$
$$\le (1+\varepsilon) W_p^d(Q_i(\cdot / \tilde{x}^{i-1}), P_i(\cdot / x^{i-1}))^p$$

for all $\tilde{x}^{i-1}, x^{i-1}$ in $E^{i-1}$. Obviously, $\tilde{X}^n, X^n$ are of law $\mathbb{Q}, \mathbb{P}$, respectively.

By the triangle inequality for the $W_p^d$-distance,

$$\mathbb{E}(d(\tilde{X}_i, X_i)^p / \tilde{X}^{i-1} = \tilde{x}^{i-1}, X^{i-1} = x^{i-1})$$
$$\le (1+\varepsilon)[W_1^d(Q_i(\cdot / \tilde{x}^{i-1}), P_i(\cdot / \tilde{x}^{i-1})) + W_1^d(P_i(\cdot / \tilde{x}^{i-1}), P_i(\cdot / x^{i-1}))]^p.$$

Using the elementary inequality that $(x+y)^p \le a^{p-1}x^p + b^{p-1}y^p$ (for $p \ge 1$ $\forall x, y \ge 0$) where $a, b > 1$ such that $1/a + 1/b = 1$, we have by the assumptions $P_i(\cdot / x^{i-1}) \in T_p(C)$ and (C1)

$$
\begin{aligned}
\mathbb{E}(d^p(\tilde{X}_i, X_i) &/ \tilde{X}^{i-1} = \tilde{x}^{i-1}, X^{i-1} = x^{i-1}) \\
&\le (1+\varepsilon)\left(\sqrt{2CH_i(\tilde{x}^{i-1})} + \left[\sum_{j=1}^{i-1}(a_j)^p d^p(\tilde{x}_{i-j}, x_{i-j})\right]^{1/p}\right)^p \\
&\le (1+\varepsilon)\left(a^{p-1}[2CH_i(\tilde{x}^{i-1})]^{p/2} + b^{p-1}\sum_{j=1}^{i-1}(a_j)^p d^p(\tilde{x}_{i-j}, x_{i-j})\right).
\end{aligned}
$$

(2.6)

By recurrence on $i$, this entails that $\mathbb{E}d^p(X_i, \tilde{X}_i) < +\infty$ for all $i = 1, \dots, n$. Taking the average with respect to $\mathcal{L}(\tilde{X}^{i-1}, X^{i-1})$, summing on $i$ and using the concavity of the function $x \to x^{p/2}$ for $p \in [1, 2]$, we get by (2.5) and (2.6)

$$\frac{1}{n(1+\varepsilon)}\sum_{i=1}^{n}\mathbb{E}(d^p(\tilde{X}_i, X_i))$$

$$\leq a^{p-1}\left(\frac{2C}{n}\sum_{i=1}^{n}\mathbb{E}H_i(\tilde{X}^{i-1})\right)^{p/2}+\frac{b^{p-1}}{n}\sum_{i=1}^{n}\sum_{j=1}^{i-1}a_j^p\mathbb{E}\,d^p(\tilde{X}_{i-j},X_{i-j})$$

$$=a^{p-1}\left(\frac{2C}{n}H(\mathbb{Q}/\mathbb{P})\right)^{p/2}+\frac{b^{p-1}}{n}\sum_{k=1}^{n-1}\mathbb{E}\,d^p(\tilde{X}_k,X_k)\sum_{i=k+1}^{n}a_{i-k}^p.$$

Using $\sum_{j\geq 1}a_j^p=r^p$ and letting $\varepsilon\to 0+$, the above inequality gives us, when $r^p b^{p-1}<1$,

$$W_p^{d_{l_p}}(\mathbb{Q},\mathbb{P})\leq\left(\frac{a^{p-1}}{1-r^p b^{p-1}}\right)^{1/p}\sqrt{2Cn^{2/p-1}H(\mathbb{Q}/\mathbb{P})}.$$

Optimizing on $(a,b)$, we get the desired inequality. $\quad\square$

Noting that for a real function $f$ on $E^n$, $\|f\|_{\mathrm{Lip}(d_{l_1})}\leq\alpha$ if and only if for every $k=1,\dots,n$,

(2.7)        $$|f_k(x_k)-f_k(y_k)|\leq\alpha d(x_k,y_k)\qquad\forall x_k,y_k\in E,$$

where $f_k(x_k)$ is the function $f$ w.r.t. the $k$th variable while the others are fixed. Then we get by Theorem 1.1,

COROLLARY 2.6.   *Under the assumption of Theorem* 2.5 *for* $p=1$, *for any real function* $f$ *on* $E^n$ *satisfying* (2.7),

$$\mathbb{E}_{\mathbb{P}}e^{\lambda(f-\mathbb{E}_{\mathbb{P}}f)}\leq\exp\left(\frac{C\lambda^2\alpha^2 n}{2(1-r)^2}\right)\qquad\forall\lambda\in\mathbb{R}.$$

*In particular, for any* $t>0$,

$$\mathbb{P}(f>\mathbb{E}_{\mathbb{P}}f+t)\leq\exp\left(-\frac{t^2(1-r)^2}{2nC\alpha^2}\right).$$

REMARK 2.7.   The condition $P_i(\cdot/x^{i-1})\in T_p(C)$ is our starting point for tensorization of the $T_p(C)$ and it is verified for many interesting examples, such as the stochastic differential equation (SDE) (4.1) or random dynamical systems or Gibbs fields. Condition (C1), meaning that the dependence of the *present* on the past is very weak, is a crucial condition. Indeed, when $d(x,y)=\mathbb{1}_{x\neq y}$, $p=1$ and $\mathbb{P}$ is Markovian, (C1) is equivalent to (1.6), and Theorem 2.5 is exactly the result of Marton mentioned in the Introduction.

REMARK 2.8.   That the constant $C_n$ for the $T_1$-inequality of $\mathbb{P}_x$ increases linearly on dimension $n$ is natural in the point of view of the Hoeffding inequality in Corollary 2.6. This is completely different from the case of the $T_2$-inequality, for which it is hoped that the $T_2$-constant remains independent of dimension $n$, as seen for the independent tensorization of $T_2(C)$ by Talagrand [18] or its extension Theorem 2.5.

REMARK 2.9. Under $P_i(\cdot / x^{i-1}) \in T_p(C)$ and (2.4) but without the contraction condition that $r^p := \sum_j (a_j)^p < 1$, we have always $\mathbb{P}_x \in T_p(C_n)$ on $E^n$ w.r.t. $d_{l_p}$ for some constant $C_n$ (but the crucial estimate of $C_n$ in Theorem 2.5 is lost). We give only the proof of this fact for $p = 1$.

Indeed, consider the nonnegative nilpotent lower triangular matrix $A = (a_{ij})$, where $a_{ij} = a_{i-j}$ if $i > j$ and 0 otherwise. For any given $\delta \in (0, 1)$, there is always a (positive) vector $z = (z_1, \ldots, z_n)$ such that $z_i > 0$, $\sum_i z_i = 1$ and

$$(zA)_k = \sum_{i=k+1}^n z_i a_{i-k} \le \delta z_k \qquad \forall k = 1, \ldots, n.$$

Then by (2.5) for $p = 1$, we have by Jensen's inequality,

$$
\begin{aligned}
\mathbb{E} \sum_{i=1}^n & z_i \, d(X_i, \tilde{X}_i) \\
&\le (1 + \varepsilon) \left( \sum_{i=1}^n z_i \mathbb{E} \sqrt{2 C H(\tilde{x}^{i-1})} + \sum_{i=1}^n z_i \sum_{j=1}^{i-1} a_j \mathbb{E} \, d(X_{i-j}, \tilde{X}_{i-j}) \right) \\
&\le (1 + \varepsilon) \left( \sqrt{\sum_{i=1}^n z_i 2 C \mathbb{E} H(\tilde{x}^{i-1})} + \sum_{k=1}^{n-1} \mathbb{E} \, d(X_k, \tilde{X}_k) \sum_{i=k+1}^n z_i a_{i-k} \right) \\
&\le (1 + \varepsilon) \left( \sqrt{2 C \max_i z_i H(\mathbb{Q}/\mathbb{P})} + \sum_{k=1}^{n-1} \delta z_k \mathbb{E} \, d(X_k, \tilde{X}_k) \right),
\end{aligned}
$$

where it follows that

$$W_1^{d_{l_1}}(\mathbb{Q}, \mathbb{P}) \le \frac{1}{(1 - \delta) \min_i z_i} \sqrt{2 C \max_i z_i H(\mathbb{Q}/\mathbb{P})}.$$

When $z_i = 1/n$, the best choice of $\delta$ is $r$, and this inequality becomes Theorem 2.5.

### 2.5. $T_1(C)$ for weakly dependent sequences: McDiarmid–Rio's martingale method revisited.

The last inequality in Corollary 2.6, applied to $F(X_1, \ldots, X_n) = \sum_{k=1}^n f(X_k)$ and the trivial metric $d$, where $(X_k)$ are independent and $\|f(X_k)\|_\infty \le \alpha$, becomes exactly the sharp Hoeffding inequality (see [13]). But when it is applied to $F(X_1, \ldots, X_n) = f(X_n)$, it does not furnish the good order of $n$ for $n$ large. As this last question is important for the $T_1(C)$ of the the invariant measure, we give now a very simple proof of the following:

PROPOSITION 2.10. *Let $(E, d)$ be a Polish space. Let $P(x, dy)$ be a Markov kernel on $E$ such that:*

(a) *$P(x, \cdot) \in T_1(C)$ for every $x \in E$;*

(b) $W_1^d(P(x, \cdot); P(\tilde{x}, \cdot)) \leq r d(x, \tilde{x})$, *for every* $x, \tilde{x}$ *in* $E$ *and some* $r < 1$.

*Then there is a unique invariant probability measure* $\mu$ *of* $P$ *and it satisfies* $T_1(C_\infty)$ *as well as* $P^n(x, \cdot) \ \forall n \geq 1$, *where* $C_\infty = C(1 - r^2)^{-1}$.

PROOF.    When $(E, d)$ is Polish, the space $M_1^p(E)$ of probability measures $\nu$ on $E$ such that $\int d(x, x_0)^p \, d\nu(x) < +\infty$, equipped with the Wasserstein metric $W_p(\cdot, \cdot)$ is a metric complete separable space (see [19]). Since $\nu \in M_1^1(E) \Rightarrow \nu P \in M_1^1(E)$ by (a) and, condition (b) implies (in fact, equivalent to)

$$W_1(\nu_1 P, \nu_2 P) \leq r W_1(\nu_1, \nu_2) \qquad \forall \nu_1, \nu_2 \in M_1^1(E),$$

hence, by the fixed point theorem, there is one and only one $P$-invariant measure $\mu \in M_1^p(E)$, and $P^n(x, \cdot) \to \mu$ in the metric $W_1$ for any initial point $x \in E$. The last point shows also that $\mu$ is the unique invariant probability measure of $P$ [without the restriction that $\mu \in M_1^1(E)$].

Since

$$W_1^d(\nu, \mu) = \sup_{f : \|f\|_{\mathrm{Lip}} \leq 1} \left| \int_E f \, d\nu - \int_E f \, d\mu \right|,$$

condition (b) is also equivalent to

$$\|Pf\|_{\mathrm{Lip}} \leq r \|f\|_{\mathrm{Lip}} \qquad \forall f.$$

Thus, $\|P^N f\|_{\mathrm{Lip}} \leq r^N \|f\|_{\mathrm{Lip}}$ for all $N \geq 1$. Now given a Lipschitzian function $f$, we have by (a) and Bobkov–Götze's Theorem 1.1,

$$
\begin{aligned}
P^n(e^f) &\leq P^{n-1}\left[ \exp\left( Pf + \frac{C\|f\|_{\mathrm{Lip}}^2}{2} \right) \right] \\
&\leq P^{n-2}\left[ \exp\left( P^2 f + \frac{C\|f\|_{\mathrm{Lip}}^2}{2} + \frac{C\|Pf\|_{\mathrm{Lip}}^2}{2} \right) \right] \\
&\leq \cdots \\
&\leq \exp\left( P^n f + \frac{C\|f\|_{\mathrm{Lip}}^2}{2} + \frac{C\|Pf\|_{\mathrm{Lip}}^2}{2} + \cdots + \frac{C\|P^{n-1}f\|_{\mathrm{Lip}}^2}{2} \right) \\
&\leq \exp\left( P^n f + \frac{C\|f\|_{\mathrm{Lip}}^2}{2(1 - r^2)} \right).
\end{aligned}
$$

In other words, for every $x \in E$, $P^n(x, \cdot) \in T_1(C_\infty)$, where $C_\infty$ is given in the proposition. Letting $n \to \infty$, we obtain the desired result for $\mu$ by Lemma 2.2. □

We now use the martingale method of McDiarmid [14] (in the independent case) and Rio [16] (in the uniform mixing case) for extending the argument above to the process-level law $\mathbb{P}$.

THEOREM 2.11. *Let $\mathbb{P}$ be a probability measure on $E^n$ satisfying $P_i(\cdot/x^{i-1}) \in T_1(C)$ ($\forall i, x^{i-1}$) in Theorem 2.5. Assume instead of* (C1) *that*

(C1′) *there is some constant $S > 0$ such that for all real bounded Lipschitzian function $f(x_{k+1}, \ldots, x_n)$ with $\|f\|_{\mathrm{Lip}(d_{l_1})} \leq 1$, for all $x \in E^n$, $y_k \in E$,*

$$|\mathbb{E}_{\mathbb{P}}(f(X_{k+1}, \ldots, X_n)/X^k = x^k) - \mathbb{E}_{\mathbb{P}}(f(X_{k+1}, \ldots, X_n)/X^k = (x^{k-1}, y_k))|$$
$$\leq Sd(x_k, y_k).$$

*Then for all function, $F$ on $E^n$ satisfying* (2.7),

$$(2.8) \qquad \mathbb{E}_{\mathbb{P}} e^{\lambda(F - \mathbb{E}_{\mathbb{P}} F)} \leq \exp\left( \frac{C\lambda^2 (1+S)^2 \alpha^2 n}{2} \right) \qquad \forall \lambda \in \mathbb{R}.$$

*Equivalently, $\mathbb{P} \in T_1(C_n)$ on $(E^n, d_{l_1})$ with*

$$C_n = nC(1+S)^2.$$

PROOF. We may assume without loss of generality that $\alpha = 1$. Let $(M_k = \mathbb{E}_{\mathbb{P}}(F/X^k))_{k \geq 0}$, where $M_0 = \mathbb{E}_{\mathbb{P}} F$. It is a martingale. It is enough to show that for each $k$,

$$\mathbb{E}_{\mathbb{P}}(e^{\lambda(M_k - M_{k-1})}/X^{k-1}) \leq \exp\left( \frac{C\lambda^2(1+S)^2}{2} \right).$$

To this end, note at first by $P_i(\cdot/x^{i-1}) \in T_1(C)$ and Theorem 1.1,

$$\mathbb{E}_{\mathbb{P}}(e^{\lambda(M_k - M_{k-1})}/X^{k-1}) \leq \exp\left( \frac{C\lambda^2 b_k^2}{2} \right),$$

where

$$b_k := \sup_{x,y} \frac{|M_k(x^k) - M_k(x^{k-1}, y_k)|}{d(x_k, y_k)}.$$

But $M_k(x^k) = \int F(x^k, x_{k+1}, \ldots, x_n) \mathbb{P}(dx_{k+1}, \ldots, dx_n/x^k)$, writing $x_{k+1}^n = (x_{k+1}, \ldots, x_n)$ we have

$$|M_k(x^k) - M_k(x^{k-1}, y_k)|$$
$$\leq \left| \int (F(x^k, x_{k+1}^n) - F(x^{k-1}, y_k, x_{k+1}^n)) \mathbb{P}(dx_{k+1}^n/x^k) \right|$$
$$\quad + \left| \int F(x^{k-1}, y_k, x_{k+1}^n)(\mathbb{P}(dx_{k+1}^n/x^k) - \mathbb{P}(dx_{k+1}^n/x^{k-1}, y_k)) \right|$$
$$\leq d(x_k, y_k) + Sd(x_k, y_k).$$

Hence, $b_k \leq (1 + S)$, the desired result.    □

REMARK 2.12.   When $d(x, y) = \mathbb{1}_{x \neq y}$, $P_i(\cdot / x^{i-1}) \in T_1(1/4)$, and this result is essentially due to Rio [16]. Using a different condition than (C1′), he essentially proved that the constant $S$ in condition (C1′) verifies $S \leq 2\sum_{j=1}^{\infty} \phi_j$, where $\phi_j$ is the uniform mixing coefficient of the sequence $(X_n)$. Our proof above is, in fact, inspired by his work.

REMARK 2.13.   If the condition (C1) is viewed as a *backward* type, then (C1′) may be seen as a *forward* type. Indeed (C1′) is equivalent to

$$W_1^{d_{l_1}}(\mathbb{P}(dx_{k+1}^n / x_k, x^{k-1}), \mathbb{P}(dx_{k+1}^n / y_k, x^{k-1})) \leq S d(x_k, y_k).$$

It means intuitively that the present does not influence a lot the future of the process $\mathbb{P}$. In concrete situations (C1′) is often weaker than (C1) with $p = 1$. For example, let $(\mathbb{P}_x)$ be a uniformly ergodic (Doeblin recurrent, say) Markov chain with transition $P(x, dy)$ in the sense that $r_n := \sup_{x \in E} \|P^n(x, \cdot) - \mu\|_{\mathrm{TV}} \to 0$. As $2\phi_n \leq r_n$, we have by Rio's estimate above,

$$S \leq \sum_{n=1}^{\infty} \sup_{x \in E} \|P^n(x, \cdot) - \mu\|_{\mathrm{TV}},$$

which is finite. But Marton's condition (1.6) or (C1) means $(1/2) \sup_{x \in E} \|P^n(x, \cdot) - \mu\|_{\mathrm{TV}} \leq r^n$ for all $n \geq 1$. See also Example 3.3.

It would be very interesting to generalize Theorem 2.11 to $T_2(C)$.

## 3. Application: study of $T_1(C)$ and $T_2(C)$ for random dynamical systems.

### 3.1. $T_1(C)$.   
Let $E$ be a complete connected Riemannian manifold equipped with the Riemannian metric $d$. Consider now the nonlinear random perturbed dynamical system valued in $E$,

$$(3.1) \quad X_0(x) := x \in E, \qquad X_{n+1}(x) = F(X_n(x), W_{n+1}), \qquad n \geq 0,$$

where the noise $(W_n)_{n \geq 0}$ is a sequence of i.i.d. r.v. valued in some measurable space $(G, \mathcal{G})$, defined on some probability space $(\Omega, \mathcal{F}, \mathbb{P})$, and $F(x, w): E \times G \to E$ is measurable. Denote by $P(x, dy)$ the law of $F(x, W_1)$, and the following:

PROPOSITION 3.1.   *Assume that there exists $\delta > 0$ such that*

$$(3.2) \qquad\qquad \sup_{x \in E} \mathbb{E}(e^{\delta d(F(x, W_1), F(x, W_2))^2}) < +\infty.$$

*If there exists $0 \leq r < 1$ such that*

$$(3.3) \qquad \mathbb{E}(d(F(x, W_1), F(\tilde{x}, W_1))) \leq r d(x, \tilde{x}) \qquad \forall x, \tilde{x} \in E,$$

*or more generally for some constant $S \geq 0$,*

$$(3.4) \qquad \sum_{n=1}^{\infty} \mathbb{E}(d(X_n(x), X_n(\tilde{x}))) \leq S\, d(x, \tilde{x}) \qquad \forall x, \tilde{x} \in E,$$

*then there is some constant $C > 0$ such that for any $n \geq 1$, for every probability measure $\mathbb{Q}^n$ on $E^n$,*

$$W_1^{d_{l_1}}(\mathbb{Q}^n, \mathbb{P}_x^n) \leq \sqrt{CnH(\mathbb{Q}^n/\mathbb{P}_x^n)},$$

*where $\mathbb{P}_x^n$ is the law of $(X_k(x))_{1 \leq k \leq n}$ on $E^n$.*

PROOF. By Theorem 2.3, condition (3.2) is equivalent to "$P(x, \cdot) \in T_1(C) \; \forall x \in E$." Notice that (3.3) is equivalent to (C1) (with $p = 1$) in Theorem 2.5, and (3.4) implies trivially (C1') with the same constant $S$ in Theorem 2.11. Hence, this proposition follows from Theorems 2.5 and 2.11. □

REMARK 3.2. If the largest Lyapunov exponent in $L^1$ given by

$$\lambda_{\max}(L^1) := \lim_{n \to \infty} \left( \sup_{x \neq \tilde{x}} \frac{\mathbb{E}d(X_n(x), X_n(\tilde{x}))}{d(x, \tilde{x})} \right)^{1/n}$$

is strictly smaller than 1, then condition (3.4) is verified.

EXAMPLE 3.3 (ARMA model). To see the difference between (C1) in Theorem 2.5 and (C1') in Theorem 2.11, let us consider the ARMA model

$$X_0(x) = x, \qquad X_{n+1}(x) = AX_n(x) + W_{n+1}$$

in $E = \mathbb{R}^d$, where $A \in \mathcal{M}_{d \times d}$ (the space of $d \times d$ matrices) and $(W_n)$ is a sequence of i.i.d. r.v. with values in $G = \mathbb{R}^d$. This model is a particular case of the general model above with $F(x, w) = Ax + w$. Condition (C1), equivalent to (3.3), means that $r = \|A\| := \sup\{|Ax|; |x| \leq 1\} < 1$, however, (C1') for this linear model is equivalent to

$$r_{\mathrm{sp}}(A) := \max\{|\lambda|; \lambda \text{ is an eigenvalue in } \mathbb{C} \text{ of } A\} = \lambda_{\max}(L^1) < 1,$$

which is much weaker. This last condition is a well-known sharp sufficient condition for the ergodicity of this linear ARMA model $(X_n)$.

REMARK 3.4. For this model, the known results mentioned in the Introduction cannot be applied, for the uniform mixing condition is, in general, not satisfied when $E$ is noncompact. For example, the ARMA model with $A \neq 0$ and $W_1$ unbounded is never uniformly mixing. See [22].

3.2. $T_2(C)$.  Consider a particular case of the preceding model

$$(3.5) \qquad X_0(x) = x, \qquad X_{n+1}(x) = f(X_n(x)) + \sigma(X_n(x))W_{n+1},$$

(the discrete time SDE), that is, $F(x, w) = f(x) + \sigma(x)w$, where $E = \mathbb{R}^d$, $G = \mathbb{R}^n$, $f : \mathbb{R}^d \to \mathbb{R}^d$, $\sigma : \mathbb{R}^d \to \mathcal{M}_{d \times n}$ (the space of $d \times n$ matrices) and the noise $(W_n)_{n \in \mathbb{Z}}$ is a sequence of i.i.d. r.v. with values in $\mathbb{R}^n$ such that $\mathbb{E}W_1 = 0$. Assume that:

(i) $\mathbb{P}_W := \mathbb{P}(W_1 \in \cdot) \in T_2(C)$ on $\mathbb{R}^n$ w.r.t. the Euclidean metric;
(ii) $|\sigma(x)w| \le K|w| \ \forall (x, w) \in \mathbb{R}^d \times \mathbb{R}^n$;
(iii) for some $r \in [0, 1)$,

$$(3.6) \qquad \sqrt{|f(x) - f(\tilde{x})|^2 + \mathbb{E}|(\sigma(x) - \sigma(\tilde{x}))W_1|^2} \le r|x - \tilde{x}| \qquad \forall x, \tilde{x} \in \mathbb{R}^d.$$

Notice that conditions (i) and (ii) imply that $P(x, \cdot) \in T_2(CK^2)$ for all $x \in \mathbb{R}^d$, by Lemma 2.1; and condition (iii) implies (C1) with the same $r$ for $p = 2$. Hence, by Theorem 2.5, $\mathbb{P}_x^n \in T_2(CK^2/(1-r)^2)$. That yields, by Bobkov, Gentil and Ledoux [1], the following:

COROLLARY 3.5. *For the model* (3.5) *above assume conditions* (i)–(iii). *Then* $\mathbb{P}_x^n \in T_2(CK^2/(1-r)^2)$ *and for any measurable function* $F(x_1, \dots, x_n) \in L^1((\mathbb{R}^d)^n, \mathbb{P}_x^n)$,

$$\mathbb{E}\exp(\rho Q F(X_1(x), \dots, X_n(x))) \le \exp(\rho\mathbb{E}F(X_1(x), \dots, X_n(x))),$$

*where*

$$\rho := \frac{(1-r)^2}{CK^2}, \qquad QF(x_1, \dots, x_n) := \inf_{y \in (\mathbb{R}^d)^n}\left(F(x+y) + \tfrac{1}{2}\sum_{k=1}^n |y_k|^2\right).$$

As noted in [1], several estimates of Laplace integrals are the consequence of the functional inequality version of the $T_2(C)$ above. For instance, Corollary 6.1 in [1] says that for any convex function $F$ on $(\mathbb{R}^d)^n$,

$$\mathbb{E}_{\mathbb{P}_x^n}\exp\left(\rho\left[F - \tfrac{1}{2}\sum_{k=1}^n (\partial_k F)^2\right]\right) \le \exp(\rho\mathbb{E}_{\mathbb{P}_x^n}F).$$

REMARK 3.6. Consider the Lyapunov exponent in $L^2$,

$$\lambda_{\max}(L^2) := \lim_{n \to \infty}\left(\sup_{x \ne \tilde{x}}\frac{\mathbb{E}d(X_n(x), X_n(\tilde{x}))^2}{d(x, \tilde{x})^2}\right)^{1/n}.$$

Obviously, (3.6) implies $\lambda_{\max}(L^2) < 1$. It is then natural to ask whether $P(x, \cdot) \in T_2(C) \ \forall x$ plus $\lambda_{\max}(L^2) < 1$ do imply "$\mathbb{P}_x^n \in T_2(K)$" for some constant $K$ independent of $n$ (for which we have no answer unlike for $T_1$). Notice that for the ARMA model, $\lambda_{\max}(L^2) = \lambda_{\max}(L^1) = r_{\mathrm{sp}}(A)$.

**4. Application: study of $T_1(C)$ for paths of SDEs.**   Let us give here an application of Theorem 2.3 to SDE. Consider the SDE in $\mathbb{R}^d$,

$$(4.1) \qquad dX_t = \sigma(X_t)\,dB_t + b(X_t)\,dt, \ X_0 = x \in \mathbb{R}^d,$$

where $\sigma : \mathbb{R}^d \to \mathcal{M}_{d \times n}$, $b : \mathbb{R}^d \to \mathbb{R}^d$ and $(B_t)$ is the standard Brownian motion valued in $\mathbb{R}^n$ defined on some well filtered probability space $(\Omega, \mathcal{F}, (\mathcal{F}_t), \mathbb{P})$.

Assume that $\sigma, b$ are locally Lipschitzian and for all $x, y \in \mathbb{R}^d$,

$$(4.2) \qquad \sup_{x \in \mathbb{R}^d} \|\sigma(x)\|_{\mathrm{HS}} \leq A, \qquad \langle y - x, b(y) - b(x) \rangle \leq B(1 + |y - x|^2),$$

where $\|\sigma\|_{\mathrm{HS}} := \sqrt{\mathrm{tr}\,\sigma\sigma^t}$ is the Hilbert–Schmidt norm, $\langle x, y \rangle$ is the Euclidean inner product and $|x| := \sqrt{\langle x, x \rangle}$. It has a unique nonexplosive solution denoted by $(X_t(x))$ whose law on the space $C(\mathbb{R}^+, \mathbb{R}^d)$ of $\mathbb{R}^d$-valued continuous functions on $\mathbb{R}^+$ will be denoted by $\mathbb{P}_x$.

COROLLARY 4.1.   *Assume the conditions above. For each $T > 0$, there exists some constant $C = C(T, A, B)$ independent of initial point $x$ such that $\mathbb{P}_x$ satisfies the $T_1(C)$ for every $x \in \mathbb{R}^d$, on the space $C([0, T], \mathbb{R}^d)$ of $\mathbb{R}^d$-valued continuous functions on $[0, T]$ equipped with the metric*

$$d_T(\gamma_1, \gamma_2) := \sup_{t \in [0, T]} |\gamma_1(t) - \gamma_2(t)|.$$

PROOF.   Let $(B_t), (\tilde{B}_t)$ be two independent Brownian motions defined on some filtered probability $(\Omega, \mathcal{F}, (\mathcal{F}_t), \mathbb{P})$ and $X_t(x), \tilde{X}_t(x)$ strong solutions of (4.1), respectively, driven by $(B_t), (\tilde{B}_t)$. Put

$$\hat{X}_t := X_t(x) - \tilde{X}_t(x), \qquad \hat{b}_t := b(X_t(x)) - b(\tilde{X}_t(x))$$

$$a(\cdot) := \sigma\sigma^t(\cdot), \qquad\qquad \bar{a}_t := a(X_t(x)) + a(\tilde{X}_t(x))$$

$$L_t := \int_0^t \sigma(X_t(x))\,dB_t - \int_0^t \sigma(\tilde{X}_t(x))\,d\tilde{B}_t.$$

Then

$$\hat{X}_t = L_t + \int_0^t \hat{b}_s\,ds.$$

By Theorem 2.3, it is enough to show that there exists some positive constant $\delta = \delta(T, A, B)$ such that

$$(4.3) \qquad \mathbb{E}\exp\left(\delta \sup_{0 \leq t \leq T} |\hat{X}_t|^2\right) < +\infty.$$

Let $f(x) := h(|x|)$, where $h \in C^\infty(\mathbb{R})$ is pair and such that $h(r) = r$ for $r \geq 4$ and

$$h(r) \geq r, \qquad 0 \leq h'(r) \leq 1 \wedge r, \qquad 0 \leq h''(r) \leq 1 \qquad \forall r \in [0, 4].$$

Consider $Y_t := (1 + f(\hat{X}_t))e^{-\beta t}$, where $\beta > 0$ is a constant to be determined later. By Ito's formula,

$$dY_t = e^{-\beta t}\left(\frac{1}{2}\sum_{i,j=1}^d \bar{a}_t^{ij}\partial_i\partial_j f(\hat{X}_t) + \langle \nabla f(\hat{X}_t), \hat{b}_t\rangle\right) dt - \beta Y_t\, dt + dM_t$$

$$= e^{-\beta t}\left(\frac{1}{2}h''(|\hat{X}_t|)\frac{\langle \hat{X}_t, \bar{a}_t\hat{X}_t\rangle}{|\hat{X}_t|^2} + \frac{1}{2}h'(|\hat{X}_t|)\left(\frac{\operatorname{tr}\bar{a}_t}{|\hat{X}_t|} - \frac{\langle \hat{X}_t, \bar{a}_t\hat{X}_t\rangle}{|\hat{X}_t|^3}\right)\right.$$

$$\left. + \frac{h'(|\hat{X}_t|)}{|\hat{X}_t|}\langle \hat{X}_t, \hat{b}_t\rangle - \beta(1 + h(|\hat{X}_t|))\right) dt + dM_t,$$

where $(M_t)$ is a local martingale $(M_t)$ with $M_0 = 0$, whose quadratic variational process $[M]$ is given by

$$[M]_t = \int_0^t e^{-2\beta s}\langle \nabla f(\hat{X}_s), \bar{a}_s \nabla f(\hat{X}_s)\rangle\, ds \le 2A^2\int_0^t e^{-2\beta s}\, ds \le \frac{A^2}{\beta}.$$

Using our condition (4.2), we see that $Y_t \le 1 + h(0) + M_t$ once if

$$\beta > \max\{0, 2A^2 + B\}.$$

Fix such a $\beta$. For any $\lambda > 0$, using the exponential martingale,

$$\exp\left(\lambda M_t - \frac{\lambda^2}{2}[M]_t\right),$$

(Novikov's condition is satisfied) and Doob's maximal inequality [applied to the positive submartingale $\exp(\lambda M_t/2)$], we have

$$\mathbb{E}e^{\lambda(\sup_{t\le T} Y_t - 1 - h(0))} \le \mathbb{E}\sup_{t\le T} e^{\lambda M_t} \le 4(\mathbb{E}e^{\lambda M_T})^2 \le 4\exp\left(\frac{\lambda^2 A^2}{\beta}\right).$$

Hence, by Chebychev's inequality and an optimization of $\lambda$, we get

$$\mathbb{P}\left(\sup_{t\le T} Y_t > 1 + h(0) + r\right) \le 4\exp\left(-\frac{\beta r^2}{4A^2}\right) \qquad \forall r > 0.$$

Consequently,

$$\mathbb{E}\exp\left(a\sup_{t\le T} Y_t^2\right) < +\infty, \qquad \text{if } 0 < a < \frac{\beta}{4A^2}.$$

Hence, (4.3) is true for all $\delta \in (0, e^{-\beta T}\frac{\beta}{4A^2})$, where $\beta > \max\{0, 2A^2 + B\}$. $\quad\square$

REMARK 4.2. If $b \in C^2$ verifies for some constant $B$,

$$(4.4) \qquad \nabla^s b := (\tfrac{1}{2}(\partial_i b^j + \partial_j b^i))_{1\le i,j\le d} \le BI_d$$

in the order of nonnegative definiteness where $I_d$ is the identity matrix, then $\langle y - x, b(y) - b(x)\rangle \le B|x - y|^2$ and the condition on $b$ in (4.2) is satisfied.

REMARK 4.3.   Assume $\|\nabla b\| \leq K$, $n = d$ and $\sigma(x) = \sigma = I_d$. Capitaine, Hsu and Ledoux [3] yields the log-Sobolev inequality below:

$$\int_{C([0,T],\mathbb{R}^d)} F^2 \log \frac{F^2}{\mathbb{E}_{\mathbb{P}_x} F^2} \, d\mathbb{P}_x \leq 2e^{KT} \int_{C([0,T],\mathbb{R}^d)} |DF|_H^2 \, d\mathbb{P}_x,$$

where $DF$ be the Malliavin gradient and

$$H := \left\{ \gamma(\cdot) := \int_0^\cdot h(s) \, ds; \|\gamma\|_H^2 = \int_0^T |h(s)|^2 \, ds < +\infty \right\}$$

(the Cameron–Martin space). As the result of Otto and Villani [15] suggests that the log-Sobolev inequality implies the $T_2(C)$ inequality (that is proved on the smooth Riemannian manifold), we should have $\mathbb{P}_x \in T_2(C)$ on $C([0,T])$ w.r.t. the following pseudo-metric,

$$d_H(\gamma_1, \gamma_2) := \begin{cases} \|\gamma_1 - \gamma_2\|_H, & \text{if } \gamma_1 - \gamma_2 \in H, \\ +\infty, & \text{otherwise.} \end{cases}$$

This last pseudo metric is much larger than $d_T$ used in the Corollary above. We shall give a simple proof of this last $T_2(C)$ inequality in Section 5.

Notice that as $d_H$ above is only a pseudo-metric and $\|X.\|_H = +\infty$, a.s., Theorem 1.1 cannot be applied for $T_1(C)$ associated with $d_H$ (since its sufficient part is no longer valid) and Theorem 2.3 (whose proof is based on Theorem 1.1) is no longer true w.r.t. $d_H$.

REMARK 4.4.   Without essential change of proof, the same result holds if the locally Lipschitzian condition of $\sigma, b$ is replaced by the well posedness of the martingale problem associated with $(\sigma\sigma^t, b)$, in the sense of Stroock–Varadhan.

REMARK 4.5.   If the condition on the drift $b$ in (4.2) is substituted by $\langle x, b(x) \rangle \leq B(1 + |x|^2) \ \forall x \in \mathbb{R}^d$, then with the same proof as above, we can prove that $\mathbb{E} \exp(\delta \sup_{t \in [0,T]} |X_t(x)|^2) < +\infty$ for some $\delta > 0$ depending on initial point. Hence, $\mathbb{P}_x$ satisfies the $T_1$-inequality with a constant $C = C_x$ depending on $x$.

Note the following drawback of the previous corollary: the constant $C$ in the $T_1$ inequality obtained through Theorem 2.3 via inequality (2.2) is of order $e^{\beta T}$ which is not natural in regard of the results obtained via weakly dependent sequences. We now show how Theorem 2.5 enables us to get the correct order.

We know from Corollary 4.1 that the law of $(X_t(x))_{t \in [0,1]}$ satisfies the $T_1$-inequality with a constant $C$ independent of $x$. In other words, the transition kernel of the Markov chain $Y_n := X_{[n,n+1]}$ valued in $C([0,1], \mathbb{R}^d)$ satisfies $T_1(C)$. Let us check (C1') below.

Given two different initial points $x, \tilde{x}$, let

$$\hat{X}_t := X_t(x) - X_t(\tilde{x}),$$

$$\hat{\sigma}_t = \sigma(X_t(x)) - \sigma(X_t(\tilde{x})), \qquad \hat{b}_t = b(X_t(x)) - b(X_t(\tilde{x})).$$

By Ito's formula,

$$|\hat{X}_t|^2 = |x - \tilde{x}|^2 + \int_0^t (\operatorname{tr}(\hat{\sigma}_s \hat{\sigma}_s^t) + 2\langle \hat{X}_s, \hat{b}_s \rangle) \, ds + M_t,$$

where $(M_t)$ is a local martingale with $M_0 = 0$, whose quadratic variational process is given by

$$[M]_t = 4 \int_0^t \langle \hat{X}_s, (\hat{\sigma}_s \hat{\sigma}_s^t) \hat{X}_s \rangle \, ds.$$

Let $\hat{\tau}_n := \inf\{t \geq 0; |\hat{X}_t| \vee [M]_t = n\}$. If there is $\delta > 0$ such that

$$(4.5) \qquad \frac{1}{2} \operatorname{tr}[(\sigma(x) - \sigma(\tilde{x}))(\sigma(x) - \sigma(\tilde{x}))^t] + \langle x - \tilde{x}, b(x) - b(\tilde{x}) \rangle$$

$$\leq -\delta |x - \tilde{x}|^2 \qquad \forall x, \tilde{x} \in \mathbb{R}^d,$$

then

$$\mathbb{E}|\hat{X}_{t \wedge \hat{\tau}_n}|^2 \leq |x - \tilde{x}|^2 - 2\delta \int_0^t \mathbb{E}|\hat{X}_{s \wedge \hat{\tau}_n}|^2 \, ds.$$

This entails by Gronwall's inequality and Fatou's lemma,

$$(4.6) \qquad \mathbb{E}|X_t(x) - X_t(\tilde{x})|^2 = \mathbb{E}|\hat{X}_t|^2 \leq |x - \tilde{x}|^2 e^{-2\delta t} \qquad \forall t \geq 0.$$

Moreover, if $\sigma$ is globally Lipchitzian, then by Burkholder–Davis–Gundy's inequality and Gronwall's inequality, we obtain easily from the estimate above that

$$\mathbb{E} \sup_{t \leq s \leq t+1} |X_s(x) - X_s(\tilde{x})|^2 \leq K |x - \tilde{x}|^2 e^{-2\delta t}$$

for some constant $K$. Thus, the Markov chain $Y_n := X_{[n,n+1]}$ valued in $C([0,1], \mathbb{R}^d)$ satisfies (C1') too. Consequently, we obtain by Theorem 2.11, the following:

PROPOSITION 4.6. *Assume* (4.2), (4.5) *and* $\sigma$ *is globally Lipchitzian. Then there is some constant* $C > 0$ *such that for any* $n \geq 1$ *and any initial point* $x$, *the law* $\mathbb{P}_x$ *of* $(X_t(x))_{t \in [0,n]}$ *on* $C([0,n], \mathbb{R}^d)$ *satisfies the inequality* $T_1(C \cdot n)$ *w.r.t. the metric*

$$d(\gamma_1, \gamma_2) := \sum_{k=0}^{n-1} \sup_{k \leq t \leq k+1} |\gamma_1(t) - \gamma_2(t)|.$$

REMARK 4.7. Let $(P_t)$ be the semigroup of transition probability kernels of our diffusion $(X_t)$. Notice that under (4.5), we have (4.6) which entails not only the existence and uniqueness of the invariant probability measure $\mu$ of $(P_t)$, but also

$$W_2^d(P_t(x, \cdot), P_t(\tilde{x}, \cdot)) \leq e^{-\delta t} |x - \tilde{x}|,$$

which gives us the exponential convergence below:

$$W_2^d(P_t(x, \cdot), \mu) \leq e^{-\delta t} \left( \int |x - \tilde{x}|^2 \, d\mu(\tilde{x}) \right)^{1/2} \qquad \forall x \in \mathbb{R}^d, t > 0.$$

Let us present a Hoeffding type inequality for

$$F(\gamma) := \int_0^n V(\gamma(t)) \, dt,$$

where $V : \mathbb{R}^d \to \mathbb{R}$ satisfies $\|V\|_{\text{Lip}} \leq \alpha$. For such $V$, $\|F\|_{\text{Lip}} \leq \alpha$ w.r.t. the metric given in the proposition above. Hence, by Theorem 1.1, Proposition 4.6 entails

$$\mathbb{P} \left( \int_0^n [V(X_t(x)) - \mathbb{E}V(X_t(x))] \, dt > r \right) \leq \exp \left( -\frac{r^2}{2nC} \right) \qquad \forall r > 0.$$

## 5. A direct approach to $T_2(C)$ for SDEs via stochastic calculus.

5.1. *$T_2$-inequality of the Wiener measure w.r.t. the Cameron–Martin metric.* Let us extend the $T_2$-inequality of the Gaussian measure due to Talagrand to the Wiener measure $\mathbb{P}$ on $C([0, T], \mathbb{R}^d)$, by means of Girsanov formula. Given $\mathbb{Q} \ll \mathbb{P}$ such that $H(\mathbb{Q}/\mathbb{P}) < +\infty$, then under $\mathbb{Q}$, there exist a Brownian motion $(B_t)$ and a predictable process $(\beta_t)$ such that the coordinates system $(\gamma_t)$ of $C([0, T], \mathbb{R}^d)$ verifies

$$d\gamma_t = dB_t + \beta_t(\gamma) \, dt, \gamma_0 = 0.$$

Moreover, it is well known that [see the proof of (5.7) below in a much more complicated case]

$$(5.1) \qquad H(\mathbb{Q}/\mathbb{P}) = \tfrac{1}{2} \mathbb{E}^{\mathbb{Q}} \int_0^T |\beta_t|^2(\gamma) \, dt.$$

Consider the Girsanov transformation $\Phi(\gamma) := \gamma(\cdot) - \int_0^\cdot \beta_t(\gamma) \, dt$. Then the law of $(\gamma, \Phi(\gamma))$ under $\mathbb{Q}$ is a coupling of $(\mathbb{Q}, \mathbb{P})$. Hence, w.r.t. the Cameron–Martin metric $d_H$ given in Remark 4.2,

$$(5.2) \quad (W_2^{d_H}(\mathbb{Q}, \mathbb{P}))^2 \leq \mathbb{E}^{\mathbb{Q}} d_H(\gamma, \Phi(\gamma))^2 = \mathbb{E}^{\mathbb{Q}} \int_0^T |\beta_t|^2(\gamma) \, dt = 2H(\mathbb{Q}/\mathbb{P}),$$

that is, $\mathbb{P} \in T_2(1)$ on $(C([0,T],\mathbb{R}^d), d_H)$. We see now why this is sharp. Indeed, if $\beta_t$ is determinist (or, equivalently, $\mathbb{Q}$ is a Gaussian measure), we claim that

$$[W_2^{d_H}(\mathbb{Q}, \mathbb{P})]^2 = \int_0^T |\beta_t|^2 \, dt = 2H(\mathbb{Q}/\mathbb{P}).$$

This follows by the following observation:

LEMMA 5.1. *Let $X$ be a random variable valued in a Banach space $E$ and $H$ be a separable Hilbert space continuously embedded in $E$. Then for any element $h \in H$,*

$$W_2^{d_H}(\mathbb{P}_X, \mathbb{P}_{X+h}) = \|h\|_H,$$

*where $\mathbb{P}_X$ is the law of $X$, $d_H(x,y) := \|x - y\|_H$ if $x - y \in H$ and $+\infty$ otherwise.*

PROOF. At first $[W_2^{d_H}(\mathbb{P}_X, \mathbb{P}_{X+h})]^2 \le \mathbb{E}\|X - (X+h)\|_H^2 = \|h\|_H^2$. To show the inverse inequality, let $\pi$ be a probability measure on $E^2$ such that its marginal laws are, respectively, laws of $X$ and $X + h$, and $\iint \|y - x\|_H^2 \pi(dx, dy) < +\infty$. Since $y - (x+h)$ is centered in the sense that $\mathbb{E}^\pi \langle e_i, y - (x+h) \rangle_H = 0$ where $(e_i)$ is an orthonormal basis of $H$, we have by Jensen's inequality,

$$\iint \|y - x\|_H^2 \pi(dx, dy) = \iint \|h + (y - (x+h))\|_H^2 \pi(dx, dy) \ge \|h\|_H^2,$$

the desired result.  □

Considering the mapping $\Psi(\gamma) = \gamma(T)$, which verifies

$$|\Psi(\gamma_1) - \Psi(\gamma_2)| \le \sqrt{T} d_H(\gamma_1, \gamma_2),$$

we get by Lemma 2.1 and (5.2) that $\mathcal{N}(0, TI_d) \in T_2(C)$ on $\mathbb{R}^d$ w.r.t. the Euclidean metric with the sharp constant $C = T$ (the theorem of Talagrand).

REMARK 5.2. Gentil [7] proved the dual (functional) version of the $T_2$-inequality of the Wiener measure w.r.t. the Cameron–Martin metric by generalizing the approach in [1]. The proof here is completely different and seems to be simpler and direct.

REMARK 5.3. Recall the method of Talagrand for proving his $T_2(C)$ for $\mathcal{N}(0, I_d)$. At first by independent tensorization, he reduces to dimension 1. And in dimension one, he uses the optimal transportation of Fréchet putting forward $\gamma = \mathcal{N}(0,1)$ to $f \, d\gamma$, and a direct integration by parts yields miraculously his $T_2(C)$. The method here is completely different, we use the Girsanov transformation which puts $\mathbb{Q}$ back to $\mathbb{P}$ instead of an (eventual) optimal transportation putting $\mathbb{P}$ forward to $\mathbb{Q}$. The approach of Talagrand is generalized recently by Feyel and Ustunel [6] who succeed to construct the optimal transportation from $\mathbb{P}$ to $\mathbb{Q}$ on an abstract Wiener space $(W, H, \mathbb{P})$.

We learned very recently (10 monthes after our first version) from Fang that the method of Girsanov transformation here has been used by Feyel and Ustunel [5] in a less elementary manner. So the result of this paragraph is due to them.

### 5.2. $T_2$-inequality of diffusions w.r.t. the Cameron–Martin metric. 

We now generalize the preceding argument to solution of the SDE

$$dX_t = dB_t + b(X_t)\,dt, \qquad X_0 = x \in \mathbb{R}^d,$$

where $(B_t)$ is a $\mathbb{R}^d$-valued Brownian motion. We assume that $b \in C^1$ and

$$\|\nabla b\| \le K.$$

For any path $\gamma \in C([0, T], \mathbb{R}^d)$ with $\gamma(0) = 0$, let $\Phi(\gamma) = \eta$ be the solution of

$$\eta(t) = x + \gamma(t) + \int_0^t b(\eta(s))\,ds.$$

Then the solution of the SDE above is given by $X_{\boldsymbol{\cdot}} = \Phi(B_{\boldsymbol{\cdot}})$. Hence, for proving the $T_2$-inequality of $X_{\boldsymbol{\cdot}}$ w.r.t. the metric $d_H$, it is enough to show that $\Phi$ is $d_H$-Lipschitzian. To this end, consider

$$g(t) = \frac{d}{d\varepsilon}\Phi(\gamma + \varepsilon h)|_{\varepsilon=0},$$

where $h \in H$ is fixed. It satisfies

$$g(t) = h(t) + \int_0^t \nabla b(\eta(s))g(s)\,ds.$$

Its solution is given by

$$g(t) = \int_0^t J(s, t)h'(s)\,ds,$$

where $J(s, t)$ is the solution of the matrix differential equation

$$(5.3) \qquad J(s, s) = I_d, \qquad \frac{d}{dt}J(s, t) = \nabla b(\eta(t))J(s, t).$$

Since $\nabla^s b \le B I_d$ for some $B \le K$, we have $|J(s, t)y| \le e^{B(t-s)}|y| \;\; \forall y \in \mathbb{R}^d$. Consequently,

$$|g(t)| \le \int_0^t e^{B(t-s)}|h'(s)|\,ds.$$

Thus, by Cauchy–Schwarz,

$$\|g\|_H^2 \le 2\int_0^T |h'(t)|^2\,dt + 2\int_0^T |\nabla b(\eta(t))g(t)|^2\,dt$$

$$\le 2\|h\|_H^2 + 2K^2\int_0^T \left[\int_0^t e^{B(t-s)}|h'(s)|\,ds\right]^2 dt.$$

Note that

$$\int_0^T \left[\int_0^t e^{B(t-s)}|h'(s)|\,ds\right]^2 dt = \int_0^T \int_0^T |h'(u)||h'(v)|\left[\int_{u\vee v}^T e^{2Bt-(u+v)}\,dt\right] du\,dv$$

$$= \langle \Gamma|h'|, |h'|\rangle_{L^2([0,T])},$$

where

$$\Gamma(u,v) = \begin{cases} e^{-B(u+v)}\dfrac{e^{2BT} - e^{2B(u\vee v)}}{2B}, & \text{if } B \ne 0, \\ T - u\vee v, & \text{if } B = 0, \end{cases}$$

and $\Gamma f(u) := \int_0^T \Gamma(u,v)f(v)\,dv$. Let $\lambda_{\max}(\Gamma)$ be the largest eigenvalue of $\Gamma$ in $L^2([0,T])$. We have $\lambda_{\max}(\Gamma) \le \|\Gamma\|_1$, the norm of $\Gamma$ in $L^1([0,T])$. It is easy to get $\|\Gamma\|_1 \le \frac{1}{B^2}$ if $B < 0$, $\|\Gamma\|_1 \le \frac{e^{2BT}}{2B^2}$ if $B > 0$, and $\|\Gamma\|_1 = \frac{T^2}{2}$ if $B = 0$. Thus, setting

$$(5.4)\qquad \alpha^2 := \alpha^2(T,K,B) = \begin{cases} 2\left(1 + \dfrac{K^2}{B^2}\right), & \text{if } B < 0, \\ 2\left(1 + K^2\dfrac{e^{2BT}}{2B^2}\right), & \text{if } B > 0, \\ 2\left(1 + \dfrac{K^2T^2}{2}\right), & \text{if } B = 0; \end{cases}$$

we get by the estimates above that $\|g\|_H^2 \le \alpha^2\|h\|_H^2$, that is, $\|\Phi\|_{\mathrm{Lip}(d_H)} \le \alpha$. Thus, Lemma 2.1 (which remains valid for the pseudo-metric $d_H$) together with the $T_2$-inequality for the Wiener measure gives us the following:

PROPOSITION 5.4.   *Assume $\nabla^s b \le BI_d$ and $\|\nabla b\| \le K$, then for every initial point $x$, $\mathbb{P}_x \in T_2(\alpha^2)$ on $C([0,T],\mathbb{R}^d)$ w.r.t. the metric $d_H$, where $\alpha^2$ is given by (5.4).*

REMARK 5.5.   Of course, the estimate of $\|\Phi\|_{\mathrm{Lip}(d_H)} \le \alpha$ together with the log-Sobolev inequality of Gross for the Wiener measure gives us also

$$\int_{C([0,T],\mathbb{R}^d)} F^2 \log\frac{F^2}{\mathbb{E}_{\mathbb{P}_x}F^2}\,d\mathbb{P}_x \le 2\alpha^2\int_{C([0,T],\mathbb{R}^d)} |DF|_H^2\,d\mathbb{P}_x,$$

which is better than the Capitaine–Hsu–Ledoux's estimate in Remark 4.3 when $B < 0$.

It is interesting to investigate whether this proposition and the corresponding log-Sobolev inequality continue to hold in the case where $\nabla^s b \leq BI_d$ with $B \leq 0$ without condition $\|\nabla b\| \leq K$.

5.3. $T_2$-inequality of diffusions w.r.t. the $L^2$-metric. Perhaps the most elementary metric on $C([0,T],\mathbb{R}^d)$ is the following $L^2[0,T]$-metric,

$$d_2(\gamma_1,\gamma_2) := \sqrt{\int_0^T |\gamma_1(t) - \gamma_2(t)|^2\,dt}.$$

Indeed, the argument leading to the $T_2$-inequality of the Wiener measure will yield the following robust $T_2$-inequality w.r.t. the metric above:

THEOREM 5.6. *Assume that $\sigma$, $b$ are locally Lipschitzian and satisfy* (4.5) *for some $\delta > 0$, and $\|\sigma\|_\infty := \sup\{|\sigma(x)z|;\ x \in \mathbb{R}^d, |z| \leq 1\} < +\infty$. Then $\mathbb{P}_x \in T_2(C)$ on $C([0,T],\mathbb{R}^d)$ w.r.t. the $L^2$-metric $d_2$ above for all $x \in \mathbb{R}^d$ and $T > 0$, where the constant $C$ is given by*

$$C := \frac{\|\sigma\|_\infty^2}{\delta^2}.$$

*Moreover, $P_T(x,\cdot) \in T_2(\frac{\|\sigma\|_\infty^2}{2\delta})$ on $\mathbb{R}^d$, as well as the unique invariant probability measure $\mu$ of $(P_t)$.*

REMARK 5.7. The two $T_2$-inequalities in this theorem are both sharp. Indeed, let $d = 1$, $\sigma(x) = 1$, $b(x) = x/2$, that is, $(X_t)$ is the standard real Ornstein–Uhlenbeck process, whose invariant measure is $\mathcal{N}(0,1)$. By this proposition, $\mu \in T_2(C)$ with $C = \|\sigma\|_\infty^2/2\delta = 1$, which is sharp.

For the sharpness of the $T_2$-inequality for $\mathbb{P}_x$ w.r.t. $d_2$, note that any Gaussian measure $\mathcal{N}(m,\Sigma)$ on $\mathbb{R}^n$ satisfies $T_2(C)$ with the sharp constant $C$ being the largest eigenvalue $\lambda_{\max}(\Sigma)$ of the covariance matrix $\Sigma$. This can be extended easily to any Gaussian measure $\nu = \mathcal{N}(m,\Sigma)$ on any separable Hilbert space $G$, where the covariance matrix $\Sigma$ is a Hilbert–Schmidt operator on $G$. Hence, if $(X_t)_{t\geq 0}$ is a Gaussian process with paths a.s. in $L^2([0,T],dt)$, then its law $\mathbb{P}$ satisfies the $T_2(C)$ on $L^2([0,T],dt)$ with the sharp constant $C = \lambda_{\max}(\Sigma)$, the largest eigenvalue of the operator

$$\Sigma f(s) := \int_0^T \mathrm{Cov}(X_s, X_t) f(t)\,dt \qquad \forall f \in L^2([0,T],dt).$$

For the Ornstein–Uhlenbeck process law $\mathbb{P}_0$ above starting from 0, $\mathrm{Cov}(X_s, X_t) = \exp(-|t-s|/2) - \exp(-(s+t)/2)$. In that case,

$$\lambda_{\max}(\Sigma) \geq \frac{\langle \Sigma \mathbb{1}_{[0,T]}, \mathbb{1}_{[0,T]}\rangle}{T} \to 4 \qquad \text{as } T \to \infty.$$

Hence, the constant $C = \|\sigma\|^2/\delta^2 = 4$ in the $T_2$-inequality for $\mathbb{P}_0$ given by our theorem becomes sharp when $T \to +\infty$.

PROOF. We shall prove that for any $\varepsilon > 0$, for any probability measure $\mathbb{Q}$ on $C([0, T], \mathbb{R}^d)$,

$$(5.5) \qquad (W_2^{d_2}(\mathbb{Q}, \mathbb{P}_x))^2 \leq 2 \frac{(1 - e^{(\varepsilon - 2\delta)T}) \|\sigma\|_\infty^2}{\varepsilon(2\delta - \varepsilon)} H(\mathbb{Q}/\mathbb{P})$$

and for any probability measure $\nu$ on $\mathbb{R}^n$,

$$(5.6) \qquad (W_2^{d_2}(\nu, P_T(x, \cdot)))^2 \leq 2 \frac{\sup_{t \in [0, T]} e^{(\varepsilon - 2\delta)t} \|\sigma\|_\infty^2}{\varepsilon} H(\nu / P_T(x, \cdot)).$$

Choosing $\varepsilon = \delta$ in (5.5), we get the first claim in the theorem; letting $\varepsilon \uparrow 2\delta$, we get $P_T(x, \cdot) \in T_2(\frac{\|\sigma\|_\infty^2}{2\delta})$ by (5.6) and then $\mu \in T_2(\frac{\|\sigma\|_\infty^2}{2\delta})$ by Lemma 2.2 and the fact that $P_T(x, \cdot) \to \mu$ as $T \to \infty$ (see Remark 4.7).

It is enough to prove (5.5) for $\mathbb{Q} \ll \mathbb{P}_x$ and $H(\mathbb{Q}/\mathbb{P}_x) < +\infty$. We divide its proof into two steps.

*Step* 1. We do at first some preparation of stochastic calculus. Let $(\Omega, \mathcal{F}, \tilde{\mathbb{P}})$ be a complete probability space on which a $n$-dimensional Brownian motion $(B_t) = (B_t^j)_{j=1,\ldots,n}$ is defined and let $\mathcal{F}_t = \mathcal{F}_t^B = \sigma(B_s, \ s \leq t)^{\tilde{\mathbb{P}}}$ (completion by $\tilde{\mathbb{P}}$). Let $X_t(x)$ be the unique solution of (4.1) starting from $x$. Then the law of $X_.(x)$ is $\mathbb{P}_x$. Consider

$$\tilde{\mathbb{Q}} := \frac{d\mathbb{Q}}{d\mathbb{P}_x}(X_.(x)) \cdot \tilde{\mathbb{P}}, \qquad M_t := \mathbb{E}^{\tilde{\mathbb{P}}} \left( \frac{d\mathbb{Q}}{d\mathbb{P}}(X_.(x)) / \mathcal{F}_t \right) \qquad \forall t \in [0, T].$$

Remark that, as $\mathbb{Q}$ is a probability measure and the law of $X(x)$ under $\tilde{\mathbb{P}}$ is exactly $\mathbb{P}_x$, we have

$$\int_\Omega \frac{d\mathbb{Q}}{d\mathbb{P}_x}(X(x)) \, d\tilde{\mathbb{P}} = \int_{C([0, T], R^d)} \frac{d\mathbb{Q}}{d\mathbb{P}_x}(w) \, d\mathbb{P}_x(w) = \mathbb{Q}(C([0, T], \mathbb{R}^d)) = 1.$$

$(M_t)$ is a martingale can and will be chosen as a continuous martingale. Let $\tau := \inf\{t \in [0, T]; \ M_t = 0\}$ with the convention that $\inf \varnothing := T+$, where $T+$ is an artificial added element larger than $T$, but smaller than any $a > T$. Then $\tilde{\mathbb{Q}}(\tau = T+) = 1$ and

$$M_t = \mathbb{1}_{t < \tau} \exp(L_t - \tfrac{1}{2}[L]_t),$$

where $L_t := \int_0^t \frac{dM_s}{M_s} \ \forall t < \tau$. $(L_t)$, being a $\tilde{\mathbb{P}}$-local martingale on $[0, \tau)$, can be represented in the following way: there is a predictable process $(\beta_t) = (\beta_t^j)_{0 \leq t < \tau}$ such that $\int_0^t |\beta_s|^2 \, ds < +\infty$, $\tilde{\mathbb{P}}$-a.s. on $[t < \tau]$ and

$$L_t = \sum_{j=1}^n \int_0^t \beta_s^j \, dB_s^j = \int_0^t \langle \beta_s, dB_s \rangle \qquad \forall t < \tau.$$

Let $\tau_n = \inf\{t \in [0, \tau[; \; [L]_t = n\}$ with the same convention that $\inf \varnothing := T+$. It is elementary that $\tau_n \uparrow \tau$, $\tilde{\mathbb{P}}$-a.s. Hence, by martingale convergence,

$$H(\mathbb{Q}/\mathbb{P}) = H(\tilde{\mathbb{Q}}/\tilde{\mathbb{P}}) = \mathbb{E}^{\tilde{\mathbb{P}}} M_T \log M_T = \lim_{n \to \infty} \mathbb{E}^{\tilde{\mathbb{P}}} M_{T \wedge \tau_n} \log M_{T \wedge \tau_n}$$

$$= \lim_{n \to \infty} \mathbb{E}^{\tilde{\mathbb{Q}}} (L_{T \wedge \tau_n} - \tfrac{1}{2} [L]_{T \wedge \tau_n}).$$

By Girsanov's formula, $(L_{t \wedge \tau_n} - [L]_{t \wedge \tau_n})_{t \in [0,T]}$ is a $\tilde{\mathbb{Q}}$-local martingale, then a true martingale since its quadratic variation process under $\tilde{\mathbb{Q}}$, being again $([L]_{t \wedge \tau_n})$, is bounded by $n$. Consequently, $\mathbb{E}^{\tilde{\mathbb{Q}}}(L_{T \wedge \tau_n} - [L]_{T \wedge \tau_n}) = 0$. Substituting it into the preceding equality and noting that $\tilde{\mathbb{Q}}(\tau_n \uparrow \tau = T+) = 1$, we get by monotone convergence,

$$(5.7) \qquad H(\mathbb{Q}/\mathbb{P}) = \tfrac{1}{2} \lim_{n \to \infty} \mathbb{E}^{\tilde{\mathbb{Q}}}[L]_{T \wedge \tau_n} = \tfrac{1}{2} \mathbb{E}^{\tilde{\mathbb{Q}}}[L]_T = \tfrac{1}{2} \mathbb{E}^{\tilde{\mathbb{Q}}} \int_0^T |\beta_t|^2 \, dt.$$

Notice that this is an extension of (5.1).

*Step* 2.  By Girsanov's theorem,

$$\tilde{B}_t := B_t - \int_0^t \beta_s \, ds$$

is a $\tilde{\mathbb{Q}}$-local martingale with $[\tilde{B}^i, \tilde{B}^j]_t = [B^i, B^j]_t = \mathbb{1}_{i=j} t$, hence, a Brownian motion under $\tilde{\mathbb{Q}}$. Under $\tilde{\mathbb{Q}}$, $X_t = X_t(x)$ verifies

$$dX_t = \sigma(X_t) \, d\tilde{B}_t + b(X_t) \, dt + \sigma(X_t)\beta_t \, dt, \qquad X_0 = x.$$

We now consider the solution $Y_t$ (under $\tilde{\mathbb{Q}}$) of

$$dY_t = \sigma(Y_t) \, d\tilde{B}_t + b(Y_t) \, dt, \qquad Y_0 = x.$$

The law of $(Y_t)_{t \in [0,T]}$ under $\tilde{\mathbb{Q}}$ is exactly $\mathbb{P}_x$. In other words, $(X, Y)$ under $\tilde{\mathbb{Q}}$ is a coupling of $(\mathbb{Q}, \mathbb{P}_x)$.

Setting

$$\hat{X}_t := X_t - Y_t, \qquad \hat{\sigma}_t := \sigma(X_t) - \sigma(Y_t), \qquad \hat{b}_t := b(X_t) - b(Y_t),$$

we have

$$(5.8) \qquad d|\hat{X}_t|^2 = [2\langle \hat{X}_t, \hat{b}_t + \sigma(X_t)\beta_t \rangle + \mathrm{tr}(\hat{\sigma}_t \sigma_t^t)] \, dt + 2\langle \hat{X}_t, \hat{\sigma}_t \, d\tilde{B}_t \rangle.$$

Letting $\hat{\tau}_n := \inf\{t \in [0,T]; |\hat{X}_t| = n\}$, we have that for any $\varepsilon > 0$,

$$\mathbb{E}^{\tilde{\mathbb{Q}}} |\hat{X}_{t \wedge \hat{\tau}_n}|^2 \leq -2\delta \int_0^t \mathbb{E}^{\tilde{\mathbb{Q}}} |\hat{X}_{s \wedge \hat{\tau}_n}|^2 \, ds + 2 \mathbb{E}^{\tilde{\mathbb{Q}}} \int_0^{t \wedge \hat{\tau}_n} \langle \hat{X}_s, \sigma(X_s)\beta_s \rangle \, ds$$

$$\leq (\varepsilon - 2\delta) \int_0^t \mathbb{E}^{\tilde{\mathbb{Q}}} |\hat{X}_{s \wedge \hat{\tau}_n}|^2 \, ds + \frac{1}{\varepsilon} \mathbb{E}^{\tilde{\mathbb{Q}}} \int_0^t \|\sigma\|_\infty^2 |\beta_s|^2 \, ds.$$

Gronwall's lemma, together with Fatou's lemma, gives us

$$(5.9) \qquad \mathbb{E}^{\tilde{\mathbb{Q}}}|\hat{X}_t|^2 \le \frac{\|\sigma\|_\infty^2}{\varepsilon} \mathbb{E}^{\tilde{\mathbb{Q}}} \int_0^t e^{(\varepsilon-2\delta)(t-s)} |\beta_s|^2 \, ds \qquad \forall \, t > 0.$$

Thus,

$$\begin{aligned}
(W_2^{d_2}(\mathbb{Q}, \mathbb{P}_x))^2 &\le \mathbb{E}^{\tilde{\mathbb{Q}}} \int_0^T |\hat{X}_t|^2 \, dt \\
&\le \frac{\|\sigma\|_\infty^2}{\varepsilon} \mathbb{E}^{\tilde{\mathbb{Q}}} \int_0^T |\beta_s|^2 \, ds \int_s^T e^{(\varepsilon-2\delta)(t-s)} \, dt \\
&\le \frac{\|\sigma\|_\infty^2}{\varepsilon} \cdot \frac{1 - e^{(2\delta-\varepsilon)T}}{2\delta - \varepsilon} \mathbb{E}^{\tilde{\mathbb{Q}}} \int_0^T |\beta_s|^2 \, ds,
\end{aligned}$$

the desired (5.5). For (5.6), notice that by the key remark (2.1),

$$H(\nu/P_T(x,\cdot)) = \inf\{H(\mathbb{Q}|_{C([0,T],\mathbb{R}^d)}/\mathbb{P}_x|_{C([0,T],\mathbb{R}^d)}); Q_T := \mathbb{Q}(x_T \in \cdot) = \nu\}.$$

And for each such $Q$, define $\tilde{Q}$ as before, we have

$$[W_2^d(\nu, P_T(x,dy))]^2 \le E^{\tilde{Q}}|\hat{X}_T|^2$$

and conclude using (5.9). □

REMARK 5.8. After the first version was submitted, we learned from M. Ledoux the work of Wang [20] who obtained the $T_2(C)$ w.r.t. the $L^2$-metric for the elliptic diffusions with lower bounded $\Gamma_2$ condition of Bakry on a Riemannian manifold. His method consists of a continuous time tensorization of the $T_2(C)$ of the heat kernels (which is true by the log-Sobolev inequality due to Bakry). Hence, the method and the result here are very different from his: the volatility coefficient $\sigma$ could be completely degenerated in Theorem 5.6, and our proof does not rely on the log-Sobolev inequality which is unknown in our context.

REMARK 5.9. By the proof above, we see that (5.5) and (5.6) hold under (4.5) even with $\delta \le 0$, except now the $T_2$-constant goes to infinity as $T \to +\infty$.

REMARK 5.10. The local Lipschitzian condition on $\sigma, b$ in this theorem can be substituted by their continuity together with the well-posedness of the martingale problem associated with $(\sigma\sigma^t, b)$. Indeed, one can find $(\sigma^n, b^n)$ tending locally uniformly to $(\sigma, b)$, such that $(\sigma^n, b^n)$ is locally Lipschitzian, $\|\sigma^n\|_\infty \le \|\sigma\|_\infty$ and verifies condition (4.5) with the same $\delta$. Now the desired result follows from Theorem 5.6 and Lemma 2.2.

As indicated in [1], many interesting consequences can be derived from this result. For instance

COROLLARY 5.11. *Under the assumptions of Theorem* 5.6, *we have for any $T > 0$,*

(a) *for any smooth cylindrical function $F$ on $G := L^2([0, T], dt; \mathbb{R}^d) \supset C([0, T], \mathbb{R}^d)$, that is,*

$$F \in \mathcal{F}C_b^\infty := \{ f(\langle \gamma, h_1 \rangle, \dots, \langle \gamma, h_n \rangle); n \geq 1, h_i \in \tilde{H}, f \in C_b^\infty(\mathbb{R}^n) \}$$

*[where $\langle \gamma_1, \gamma_2 \rangle := \int_0^T \gamma_1(t) \gamma_2(t) \, dt$], the following Poincaré inequality holds:*

$$(5.10) \qquad \mathrm{Var}_{\mathbb{P}_x}(F) \leq \frac{\|\sigma\|_\infty^2}{\delta^2} \int_{C([0, T], \mathbb{R}^d)} \|\nabla F(\gamma)\|_G^2 \, d\mathbb{P}_x(\gamma),$$

*where $\mathrm{Var}_{\mathbb{P}_x}(F)$ is the variance of $F$ under law $\mathbb{P}_x$, and $\nabla F(\gamma) \in G$ is the gradiant of $F$ at $\gamma$.*

(b) *For any $g \in C_b^\infty(\mathbb{R}^d)$,*

$$(5.11) \qquad \mathrm{Var}_{P_T(x, \cdot)}(g) \leq \frac{\|\sigma\|_\infty^2}{2\delta} \int_{\mathbb{R}^d} |\nabla g(y)|^2 P_T(x, dy).$$

(c) (*Inequality of Tsirel'son type.*) *For any nonempty subset $K$ in $G$ such that $Z(\gamma) := \sup_{h \in K} \langle \gamma, h \rangle \in L^1(\mathbb{P}_x)$, then*

$$(5.12) \qquad \int \exp \left( \frac{\delta^2}{\|\sigma\|_\infty^2} \sup_{h \in K} \left[ \langle \gamma, h \rangle - \frac{|h|_G^2}{2} \right] \right) d\mathbb{P}_x \leq \exp \left( \frac{\delta^2}{\|\sigma\|_\infty^2} \mathbb{E}^{\mathbb{P}_x} Z \right).$$

(d) (*Inequality of Hoeffding type.*) *For any $V : \mathbb{R}^d \to \mathbb{R}$ such that $\|V\|_{\mathrm{Lip}} \leq \alpha$,*

$$\mathbb{P} \left( \frac{1}{T} \int_0^T V(X_t(x)) \, dt - \mathbb{E} \frac{1}{T} \int_0^T V(X_t(x)) \, dt > r \right)$$

$$\leq \exp \left( -\frac{Tr^2 \|\sigma\|_\infty^2}{2\alpha^2 \delta^2} \right) \qquad \forall r > 0.$$

PROOF. For part (a), for any $F(\gamma) = f(\langle \gamma, h_1 \rangle, \dots, \langle \gamma, h_n \rangle) \in \mathcal{F}C_b^\infty$, we may assume without loss of generality that $h_1, \dots, h_n$ are orthonormal. In such case,

$$\Phi : \gamma \to (\langle \gamma, h_1 \rangle, \dots, \langle \gamma, h_n \rangle), \qquad G \to \mathbb{R}^n$$

is Lipschitzian with $\|\Phi\|_{\mathrm{Lip}} \leq 1$. Hence, $\nu := \mathbb{P}_x \circ \Phi^{-1} \in T_2(\|\sigma\|_\infty^2 / \delta^2)$ on $\mathbb{R}^n$ by Lemma 2.1. Thus, the result of [1], Section 4.1 entails

$$\mathrm{Var}_{\mathbb{P}_x}(F) = \mathrm{Var}_\nu(f) \leq \frac{\|\sigma\|_\infty^2}{\delta^2} \int_{\mathbb{R}^n} |\nabla f|^2 \, d\nu$$

$$= \frac{\|\sigma\|_\infty^2}{\delta^2} \int_{C([0, T], \mathbb{R}^d)} \|\nabla F(\gamma)\|_G^2 \, d\mathbb{P}_x(\gamma).$$

Part (b) is a consequence of Theorem 5.6 by [1], Section 4.1. One can derive part (c) from Theorem 5.6 by the same argument as in the finite-dimensional case given in [1], Section 6.1. For part (d), note that $T_2(C) \Rightarrow T_1(C)$. Moreover, the function $F(\gamma) := (1/T) \int_0^T V(\gamma(t)) \, dt$ on $C([0, T], \mathbb{R}^d)$ is Lipschitzian w.r.t. the $L^2$-metric and $\|F\|_{\mathrm{Lip}} \leq \alpha / \sqrt{T}$. Hence, part (d) follows from Theorem 1.1. □

REMARK 5.12. Let us compare the $T_2(C)$-inequality on $C([0, T], \mathbb{R}^d)$ w.r.t. the $L^2$-metric $d_2$ or the Cameron–Martin metric $d_H$, denoted, respectively, by $T_2(C/d_2)$, $T_2(C/d_H)$.

(a) If $\gamma_1(0) = \gamma_2(0)$, then $d_2(\gamma_1, \gamma_2) \leq \frac{2T}{\pi} d_H(\gamma_1, \gamma_2)$ by the classical Poincaré inequality. Hence, if the law $\mathbb{P}_x$ of our diffusion starting from $x$ verifies $T_2(C/d_H)$ on $C([0, T], \mathbb{R}^d)$, then $\mathbb{P}_x \in T_2(C(4T^2/\pi^2)/d_2)$ on $C([0, T], \mathbb{R}^d)$. That order $T^2$ in the last $T_2$-inequality is of correct order. For example, for the real Wiener measure $\mathbb{P}$, we see by Section 5.1 that $\mathbb{P} \in T_2(1/d_H)$ on $C([0, T], \mathbb{R}^d)$, but the largest eigenvalue $\lambda_{\max}(\Gamma)$ of the covariance function $\Gamma(s, t) = s \wedge t$ in $L^2([0, T])$ verifies

$$\lambda_{\max}(\Gamma) \geq \frac{\langle \Gamma \mathbb{1}_{[0,T]}, \mathbb{1}_{[0,T]} \rangle}{T} = \frac{T^2}{3}.$$

Thus, by the same analysis as in Remark 5.7, $\mathbb{P} \in T_2(CT^2/d_2)$ with $4/\pi^2 \geq C = \lambda_{\max}(\Gamma) \geq 1/3$.

(b) The contribution of $|\gamma_1(t) - \gamma_2(t)|$ to the $L^2$-metric is homogeneous in time $t$, but not at all to the Cameron–Martin metric $d_H$. This is the principal reason for

(b.1) The $T_2(C/d_H)$ is well adapted to the small time asymptotics of the diffusions, but not for their large time asymptotics. For instance, if $\mathbb{P}_x \in T_2(C/d_H)$, since for $Z(\gamma) = \sup_{0 \leq t \leq T} \|\gamma(t) - \gamma(0)\|$, $\|Z\|_{\mathrm{Lip}(d_H)} \leq \sqrt{T}$, then by Theorem 1.1 (its necessary part remains true for $d_H$-Lipchitzian function $F$ which is, moreover, integrable, by following the proof in [2]),

$$\mathbb{P}_x \left( \sup_{0 \leq t \leq T} |X_t(x) - x| - \mathbb{E}^x \sup_{0 \leq t \leq T} |X_t(x) - x| > r \right) \leq \exp\left( -\frac{r^2}{2CT} \right)$$

which is of the correct order when $T \to 0+$, but completely meaningless for $T$ large. See [21] for the nonadaptability of the log-Sobolev inequality w.r.t. $d_H$ for the large time asymptotics of the diffusions.

(b.2) In contrary, we have seen that the $T_2(C/d_2)$ is very well adapted for the large time asymptotics of the diffusions.

REMARK 5.13.    Theorem 5.6, together with Corollary 3.5, is our main new example for which $T_2(C)$ is true but the inequality of log-Sobolev is unknown. They are our (very partial) answer to Question 3 in the Introduction. We believe that in the situations of Theorem 5.6 and Corollary 3.5, the log-Sobolev inequality may fail without further regularity assumptions on the volatility coefficient $\sigma$.

**Acknowledgments.**    We are grateful to P. Cattiaux, S. Fang, M. Ledoux, C. Villani and F. Y. Wang for their comments on the first version. We are particularly indebted to the referee for his very numerous and conscientious remarks which led to a complete re-organization of the paper and improved the readability of the paper and for the suggestion of the simplified proofs of Lemma 2.2 and of the example at the end of Remark 2.4.

## REFERENCES

[1]  BOBKOV, S., GENTIL, I. and LEDOUX, M. (2001). Hypercontractivity of Hamilton–Jacobi equations. *J. Math. Pures Appl. (9)* **80** 669–696. MR1846020

[2]  BOBKOV, S. and GÖTZE, F. (1999). Exponential integrability and transportation cost related to logarithmic Sobolev inequalities. *J. Funct. Anal.* **163** 1–28. MR1682772

[3]  CAPITAINE, M., HSU, E. P. and LEDOUX, M. (1997). Martingale representation and a simple proof of logarithmic Sobolev inequality on path spaces. *Electron Comm. Probab.* **2** 71–81. MR1484557

[4]  DEMBO, A. (1997). Information inequalities and concentration of measure. *Ann. Probab.* **25** 927–939. MR1434131

[5]  FEYEL, D. and USTUNEL, A. S. (2002). Measure transport on Wiener space and Girsanov theorem. *C. R. Acad. Sci. Paris Sér. I Math.* **334** 1025–1028. MR1913729

[6]  FEYEL, D. and USTUNEL, A. S. (2004). The Monge–Kantorovitch problem and Monge–Ampère equation on Wiener space. *Probab. Theory Related Fields.* To appear. MR2036490

[7]  GENTIL, I. (2001). Inégalités de Sobolev logarithmiques et hypercontractivité en mécanique statistique et en E.D.P. Thèse de doctorat, Univ. Paul Sabatier Toulouse.

[8]  LEDOUX, M. (2001). *The Concentration of Measure Phenomenon.* Amer. Math. Soc., Providence, RI. MR1849347

[9]  LEDOUX, M. (2002). Concentration, transportation and functional inequalities. Preprint.

[10]  MARTON, K. (1996). Bounding $\overline{d}$-distance by information divergence: A method to prove measure concentration. *Ann. Probab.* **24** 857–866. MR1404531

[11]  MARTON, K. (1997). A measure concentration inequality for contracting Markov chains. *Geom. Funct. Anal.* **6** 556–571. MR1392329

[12]  MARTON, K. (1998). Measure concentration for a class of random processes. *Probab. Theory Related Fields* **110** 427–439. MR1616492

[13]  MASSART, P. (2003). Concentration inequalities and model selection. In *Saint-Flour Summer School.*

[14]  MCDIARMID, C. (1989). On the method of bounded differences. *Surveys of Combinatorics* (J. Siemons, ed.). *London Math. Soc. Lecture Notes Ser.* **141** 148–188. MR1036755

[15] OTTO, F. and VILLANI, C. (2000). Generalization of an inequality by Talagrand and links with the logarithmic Sobolev inequality. *J. Funct. Anal.* **173** 361–400. MR1760620

[16] RIO, E. (2000). Inégalités de Hoeffding pour les fonctions Lipschitziennes de suites dépendantes. *C. R. Acad. Sci. Paris Sér. I Math.* **330** 905–908. MR1771956

[17] SAMSON, P. M. (2000). Concentration of measure inequalities for Markov chains and $\phi$-mixing process. *Ann. Probab.* **1** 416–461. MR1756011

[18] TALAGRAND, M. (1996). Transportation cost for Gaussian and other product measures. *Geom. Funct. Anal.* **6** 587–600. MR1392331

[19] VILLANI, C. (2003). *Topics in Optimal Transportation.* Amer. Math. Soc., Providence, RI. MR1964483

[20] WANG, F. Y. (2002). Transportation cost inequalities on path spaces over Riemannian manifolds. *Illinois J. Math.* **46** 1197–1206. MR1988258

[21] WU, L. (2000). A deviation inequality for non-reversible Markov processes. *Ann. Inst. H. Poincaré Probab. Statist.* **36** 435–445. MR1785390

[22] WU, L. (2002). Essential spectral radius for Markov semigroups. I: Discrete time case. *Probab. Theory Related Fields* **128** 255–321. MR2031227

LABORATOIRE
  DE MATHÉMATIQUES APPLIQUÉES
CNRS-UMR 6620
UNIVERSITÉ BLAISE PASCAL
63177 AUBIÈRE
FRANCE
E-MAIL: djellout@math.univ-bpclermont.fr

CEREMADE
CNRS-UMR 7534
UNIVERSITÉ PARIS IX DAUPHINE
75775 PARIS
FRANCE
E-MAIL: guillin@ceremade.dauphine.fr

LABORATOIRE
  DE MATHÉMATIQUES APPLIQUÉES
CNRS-UMR 6620
UNIVERSITÉ BLAISE PASCAL
63177 AUBIÈRE
FRANCE
AND
DEPARTMENT OF MATHEMATICS
WUHAN UNIVERSITY
430072
CHINA
E-MAIL: li-ming.wu@math.univ-bpclermont.fr